\documentclass[12pt,twoside]{article}
\usepackage{latexsym,amsmath}
\usepackage{graphicx}
\def\,{\mskip 3mu} \def\>{\mskip 4mu plus 2mu minus 4mu} \def\;{\mskip 5mu plus 5mu} \def\!{\mskip-3mu}
\def\dispmuskip{\thinmuskip= 3mu plus 0mu minus 2mu \medmuskip=  4mu plus 2mu minus 2mu \thickmuskip=5mu plus 5mu minus 2mu}
\def\textmuskip{\thinmuskip= 0mu                    \medmuskip=  1mu plus 1mu minus 1mu \thickmuskip=2mu plus 3mu minus 1mu}
\textmuskip
\def\beq{\dispmuskip\begin{equation}}    \def\eeq{\end{equation}\textmuskip}
\def\beqn{\dispmuskip\begin{displaymath}}\def\eeqn{\end{displaymath}\textmuskip}
\def\bqa{\dispmuskip\begin{eqnarray}}    \def\eqa{\end{eqnarray}\textmuskip}
\def\bqan{\dispmuskip\begin{eqnarray*}}  \def\eqan{\end{eqnarray*}\textmuskip}

\newenvironment{keywords}{\centerline{\bf\small
Keywords}\begin{quote}\small}{\par\end{quote}\vskip 1ex}
\def\paradot#1{\vspace{1.3ex plus 0.5ex minus 0.5ex}\noindent{\bf{#1.}}}

\def\aidx#1{}
\def\req#1{(\ref{#1})}

\def\nq{\hspace{-1em}}

\def\odt{{\textstyle{1\over 2}}}

\def\SetR{I\!\!R}

\def\SetZ{Z\!\!\!Z}
\def\qmbox#1{{\quad\mbox{#1}\quad}}
\def\e{{\rm e}}                        

\def\E{{\bf E}}                         
\def\v{\vec}

\def\s{\sigma}
\def\v{\boldsymbol} 
\def\text#1{\mbox{\scriptsize #1}}
\def\lesssim{\mbox{\raisebox{-0.8ex}{$\stackrel{\displaystyle<}\sim$}}} 
\def\Var{\mbox{Var}}
\def\aai{{\scriptscriptstyle[]}}         

\begin{document}

\title{\vspace{-3ex}\normalsize\sc Technical Report \hfill IDSIA-14-05
\vskip 2mm\bf\Large\hrule height5pt \vskip 6mm
Bayesian Regression of \\ Piecewise Constant Functions
\vskip 6mm \hrule height2pt}
\author{{\bf Marcus Hutter}\\[3mm]
\normalsize IDSIA, Galleria 2, CH-6928\ Manno-Lugano, Switzerland\\
\normalsize marcus@idsia.ch \hspace{9ex} http://www.idsia.ch/$^{_{_\sim}}\!$marcus }
\date{10 July 2005}
\maketitle

\begin{abstract}\noindent
We derive an exact and efficient Bayesian regression algorithm for
piecewise constant functions of unknown segment number, boundary
location, and levels. It works for any noise and segment level
prior, e.g.\ Cauchy which can handle outliers. We derive simple
but good estimates for the in-segment variance. We also propose a
Bayesian regression curve as a better way of smoothing data
without blurring boundaries. The Bayesian approach also allows
straightforward determination of the evidence, break probabilities
and error estimates, useful for model selection and significance
and robustness studies. We discuss the performance on synthetic
and real-world examples. Many possible extensions will be
discussed.
\def\contentsname{\centering\normalsize Contents}
{\parskip=-2.5ex\tableofcontents}
\end{abstract}


\newpage
\section{Introduction}\label{secInt}

We consider the problem of fitting a piecewise constant function
through noisy one-dimensional data, as e.g.\ in
Figure \ref{figGLPCR}, where the segment number, boundaries and
levels are unknown. Regression with piecewise constant (PC) functions,
also known as change point detection, has many applications. For
instance, determining DNA copy numbers in cancer cells from
micro-array data, to mention just one recent.

\paradot{Bayesian piecewise constant regression (BPCR)}
We provide a full Bayesian analysis of PC-regression. For a fixed
number of segments we choose a uniform prior over all possible
segment boundary locations. Some prior on the segment levels and
data noise within each segment is assumed. Finally a prior over the
number of segments is chosen. From this we obtain the posterior
segmentation probability distribution (Section \ref{secGM}).
In practice we need summaries of this complicated distribution. A
simple maximum (MAP) approximation or mean does not work
here. The right way is to proceed in stages from determining
the most critical segment number, to the boundary location, and
finally to the then trivial segment levels. We also extract the
evidence, the boundary probability distribution, and an
interesting non-PC regression curve including error estimate
(Section \ref{secQI}).
We derive an exact polynomial-time dynamic-programming-type
algorithm for all quantities of interest (Sections \ref{secES} and
\ref{secAlg}).
Our algorithm works for any noise and level prior. We consider
more closely the Gaussian ``standard'' prior and heavy-tailed
robust-to-outliers distributions like the Cauchy, and briefly discuss
the non-parametric case (Sections \ref{secSM} and \ref{secSSD}).
Finally, some hyper-parameters like the global data average and
variability and local within-level noise have to be determined. We
introduce and discuss efficient semi-principled estimators,
thereby avoiding problematic or expensive numerical EM or
Monte-Carlo estimates (Section \ref{secHP}).
We test our method on some synthetic examples (Section
\ref{secSE}) and some real-world data sets (Section \ref{secRWE}).
The simulations show that our method handles difficult data with
high noise and outliers well.
Our basic algorithm can (easily) be modified in a variety of ways:
For discrete segment levels, segment dependent variance, piecewise
linear and non-linear regression, non-parametric noise prior, etc.
(Section \ref{secMisc}).

\paradot{Comparison to other work}
Sen and Srivastava \cite{Sen:75} developed a frequentist solution
to the problem of detecting a single (the most prominent) segment
boundary (called change or break point). Olshen et al.\
\cite{Olshen:04} generalize this method to detect pairs of break
points, which improves recognition of short segments. Both methods
are then (heuristically) used to recursively determine further
change points.
Another approach is penalized Maximum Likelihood (ML). For a fixed
number of segments, ML chooses the boundary locations that
maximize the data likelihood (minimize the mean square data
deviation). Jong et al.\ \cite{Jong:03} use a population based
algorithm as minimizer, while Picard et al.\ \cite{Picard:05} use
dynamic programming, which is structurally very close to our core
recursion, to find the exact solution in polynomial time. An
additional penalty term has to be added to the likelihood in order
to determine the correct number of segments. The most principled
penalty is the Bayesian Information Criterion
\cite{Schwarz:78,Kaass:95}. Since it can be biased towards too
simple \cite{Weakliem:99} or too complex \cite{Picard:05} models,
in practice often a heuristic penalty is used. An interesting heuristic,
based on the curvature of the log-likelihood as a function of the
number of segments, has been used in \cite{Picard:05}.
Our Bayesian regressor is a natural response to penalized ML.
Many other regressors exist; too numerous to list them all.
Another closely related work to ours is Bayesian bin density
estimation by Endres and F\"oldi\'ak \cite{Endres:05}, who also
average over all boundary locations, but in the context of density
estimation.

\paradot{Advantages of Bayesian regression}
A full Bayesian approach (when computationally feasible) has
various advantages over others: A generic advantage is that it is
more principled and hence involves fewer heuristic design choices.
This is particularly important for estimating the number of
segments. Another generic advantage is that it can be easily
embedded in a larger framework. For instance, one can decide among
competing models solely based on the (Bayesian) evidence. Finally,
Bayes often works well in practice, and provably so if the model
assumptions are valid.\footnote{Note that we are not claiming here
that BPCR works better than the other mentioned approaches. In a
certain sense Bayes is optimal if the prior is `true'. Practical
superiority likely depends on the type of application. A
comparison for micro-array data is in progress
\cite{Hutter:06genex}. The major aim of this paper is to derive an
efficient algorithm, and demonstrate the gains of BPCR beyond bare
PC-regression, e.g.\ the (predictive) regression {\em curve}
(which is better than local smoothing which wiggles more and blurs
jumps). }
We can also extract other information (nearly for free), like
probability estimates and variances for the various quantities of
interest.
Particularly interesting is the expected level (and variance) of
each data point. This leads to a regression curve, which is very
flat, i.e.\ smoothes the data, in long and clear segments, wiggles
in less clear segments, follows trends, and jumps at the segment
boundaries. It thus behaves somewhat between local smoothing
(which wiggles more and blurs jumps) and rigid PC-segmentation.

\section{The General Model}\label{secGM}

\paradot{Setup}
We are given a sequence $\v y=(y_1,...,y_n)$, e.g.\ times-series
data or measurements of some function at locations $1...n$, where
each $y_i\in\SetR$ resulted from a noisy ``measurement''., i.e.\
we assume that the $y_i$ are independently (e.g.\ Gaussian)
distributed with means $\mu_i'$ and\footnote{More generally,
$\mu'_i$ and $\s'_i$ are location and scale parameters of a
symmetric distribution.} variances $\s'_i\!\,^2$. The data
likelihood is therefore\footnote{For notational and verbal
simplicity we will not distinguish between probabilities of
discrete variables and densities of continuous variables.}
\beq\label{lh}
  \mbox{likelihood:} \qquad P(\v y|\v\mu',\v\s') \;:=\;
  \prod_{i=1}^n P(y_i|\mu'_i,\s'_i)
\eeq
The estimation of the true underlying function $f=(f_1,...,f_n)$
is called regression. We assume or model $f$ as piecewise
constant. Consider $k$ segments with segment boundaries
$0=t_0<t_1<...<t_{k-1}<t_k=n$, i.e. $f$ is constant
on $\{t_{q-1}+1,..,t_q\}$ for each $0<q\leq k$. If the noise within
each segment is the same, we have
\beq\label{pc}
  \mbox{piecewise constant:}\quad \mu'_i=\mu_q \qmbox{and} \s'_i=\s_q \qmbox{for} t_{q-1}<i\leq t_q
  \quad\forall q
\eeq
We first consider the case in which the variances of all segments
coincide, i.e.\ $\s_q=\s$ $\forall q$. Our goal is to estimate the
segment levels $\v\mu=(\mu_1,...,\mu_k)$, boundaries $\v
t=(t_0,...,t_k)$, and their number $k$. Bayesian regression
proceeds in assuming a prior for these {\em quantities of
interest}. We model the segment levels by a broad (e.g.\ Gaussian)
distribution with mean $\nu$ and variance $\rho^2$. For the
segment boundaries we take some (e.g.\ uniform) distribution among
all segmentations into $k$ segments. Finally we take some prior
(e.g.\ uniform) over the segment number $k$. So our prior
$P(\v\mu,\v t,k)$ is the product of
\beq\label{prior}
  \mbox{prior:} \qquad P(\mu_q|\nu,\rho)\,\forall q \qmbox{and}
  P(\v t|k) \qmbox{and} P(k)
\eeq
We regard the global variance $\rho^2$ and mean $\nu$ of $\v\mu$ and
the in-segment variance $\s^2$ as fixed hyper-parameters, and
notationally suppress them in the following. We will return to
their determination in Section \ref{secHP}.

\paradot{Evidence and posterior}
Given the prior and likelihood we can compute the data evidence
and posterior $P(\v y|\v\mu,\v t,k)$ by Bayes' rule:
\beqn
  \mbox{evidence:} \quad P(\v y) \;=\;
  \sum_{k,\v t} \int P(\v y|\v\mu,\v t,k)P(\v\mu,\v t,k)\, d\v\mu
\eeqn
\beqn
  \mbox{posterior:} \quad P(\v\mu,\v t,k|\v y)
  \;=\; {P(\v y|\v\mu,\v t,k)P(\v\mu,\v t,k)\over P(\v y)}
\eeqn
The posterior contains all information of interest, but is
a complex object for practical use. So we need summaries like the
maximum (MAP) or mean and variances. MAP over continuous
parameters ($\v\mu$) is problematic, since it is not
reparametrization invariant. This is particularly dangerous if MAP
is across different dimensions ($k$), since then even a linear
transformation ($\v\mu\leadsto\alpha\v\mu$) scales the posterior
(density) exponentially in $k$ (by $\alpha^k$). This
severely influences the maximum over $k$, i.e.\ the estimated
number of segments. The mean of $\v\mu$ does not have this
problem. On the other hand, the mean of $\v t$ makes only sense
for fixed (e.g.\ MAP) $k$. The most natural solution is to proceed
in stages similar to as the prior \req{prior} has been formed.

\section{Quantities of Interest}\label{secQI}

We now define estimators for all quantities of interest in stages
as suggested in Section \ref{secGM}.

\paradot{Quantities of interest}
Our first quantities are the posterior of the number of segments
and the MAP segment number
\beqn
  \mbox{\# segments:}\quad P(k|\v y) \qmbox{and} \hat k \;=\; \arg\max_k P(k|\v y)
\eeqn
Second, for each boundary $t_q$ its posterior and MAP, given the MAP estimate of $k$
\beqn
  \mbox{boundaries:}\quad P(t_q|\v y,\hat k) \qmbox{and} \hat t_q \;=\; \arg\max_{t_q} P(t_q|\v y,\hat k)
\eeqn
Different estimates of $t_q$ (e.g.\ the mean or MAP based on the
joint $\v t$ posterior) will be discussed later.
Finally we want the segment level means for the MAP segmentation
\beqn
  \mbox{segment level:}\quad P(\mu_q|\v y,\v{\hat t},\hat k) \qmbox{and}
  \hat\mu_q \;=\; \int P(\mu_q|\v y,\v{\hat t},\hat k) \mu_q d\mu_q
\eeqn
The estimate ($\hat{\v\mu},\hat{\v t},\hat k$) defines a (single) piecewise
constant (PC) function $\hat f$, which is our estimate of $f$.
A (very) different quantity is to Bayes-average over all piecewise constant
functions and to ask for the mean at location $i$ as an estimate for $f_i$.
\beqn
  \mbox{regression curve:}\quad P(\mu'_i|\v y) \qmbox{and}
  \hat\mu'_i \;=\; \int P(\mu'_i|\v y) \mu'_i d\mu'_i
\eeqn
We will see that $\v\mu'$ behaves similar to a local smoothing of
$\v y$, but without blurring true jumps. Standard deviations of
all estimates may also be reported.

\section{Specific Models}\label{secSM}

We now complete the specification of the data noise and prior.

\paradot{Segment boundaries}
We assume a uniform prior over all segmentations into $k$ segments.
Since there are $({n-1\atop k-1})$ ways of placing the $k-1$
inner boundaries (ordered and without repetition) on $(1,...,n-1)$, we have
\beq\label{ubp}
  \mbox{uniform boundary prior:}\quad
  P(\v t|k)=\textstyle({n-1\atop k-1})^{-1}
\eeq
This is the only (additional) essential assumption to be able to
derive efficient algorithms. We now discuss some (purely exemplary)
choices for the data noise and priors on $\v\mu$ and $k$.

\paradot{Gaussian model}
The standard assumption on
the noise is independent Gauss:
\beq\label{Gn}
  \mbox{Gaussian noise:}\quad
  P(y_i|\mu'_i,\s'_i) \;=\; {1\over\sqrt{2\pi}\s'_i}\;
  \mbox{\large e}^{\textstyle-{(y_t-\mu'_i)^2\over 2\s'_i\!\,^2}}
\eeq
The corresponding standard ``conjugate'' prior on the means
$\mu_q$ for each segment $q$ is also Gauss
\beq\label{Gp}
  \mbox{Gaussian prior:}\quad P(\mu_q|\nu,\rho) \;=\; {1\over\sqrt{2\pi}\rho}\;
  \mbox{\large e}^{\textstyle-{(\mu_q-\nu)^2\over 2\rho^2}}
\eeq

\paradot{Cauchy model}
The standard problem with Gauss is that it does not handle
outliers well. If we do not want to or cannot remove outliers by hand,
we have to properly model them as a prior with heavier tails.
This can be achieved by a mixture of Gaussians or by a Cauchy distribution:
\beq\label{Cn}
  \mbox{Cauchy noise:}\quad
  P(y_i|\mu'_i,\s'_i) \;=\; {1\over\pi}{\s'_i\over\s'_i\!\,^2+(y_i-\mu'_i)^2}
\eeq
Note that $\mu'_i$ and $\s'_i$ determine the location and scale of
Cauchy but are not its mean and variance (which do not exist).
The prior on the levels $\mu_q$ may as well be modeled as Cauchy:
\beq\label{Cp}
  \mbox{Cauchy prior:}\quad P(\mu_q|\nu,\rho) \;=\;
  {1\over\pi}{\rho\over\rho^2+(\mu_q-\nu)^2}
\eeq
Actually, the Gaussian noise model may well be combined with
a non-Gaussian prior and vice versa if appropriate.

\paradot{Number of segments}
Finally, consider the number of segments $k$, which is an integer
between 1 and $n$. Sure, if we have prior knowledge on the
[minimal,maximal] number of segments $[k_{min},k_{max}]$ we
could/should set $P(k)=0$ outside this interval. Otherwise, any non-extreme
choice of $P(k)$ has little influence on the final results, since
it gets swamped by the (implicit) strong (exponential) dependence
on $k$ of the likelihood. So we suggest a uniform prior
\beqn
  P(k) \;=\; {1\over k_{max}} \qmbox{for} 1\leq k\leq k_{max}
  \qmbox{and} 0 \qmbox{otherwise}
\eeqn
with $k_{max}=n$ as default (or $k_{max}<n$ discussed later).

\section{Efficient Solution}\label{secES}

\paradot{Notation}
We now derive expressions for all quantities of interest, which
need time $O(k_{max}n^2)$ and space $O(n^2)$. Throughout this
and the next section we use the following notation: %
$k$ is the total number of segments, %
$t$ some data index, %
$q$ some segment index, %
$1\leq i<h<j\leq n$ are data item indices of segment boundaries
$t_0\leq t_l<t_p<t_m\leq t_k$, %
i.e.\ $t_0=0$, $t_l=i$, $t_p=h$, $t_m=j$, $t_k=n$.
Further, $y_{ij}=(y_{i+1},...,y_j)$ is data %
with segment boundaries $t_{lm}=(t_l,...,t_m)$ %
and segment levels $\mu_{lm}=(\mu_{l+1},...,\mu_m)$.
In particular $y_{0n}=\v y$, $t_{0k}=\v t$, and $\mu_{0k}=\v\mu$.
All introduced matrices below (capital symbols with indices) will be
important in our algorithm.

\paradot{General recursion}
For $m=l+1$, $y_{ij}$ is data from a single segment with mean $\mu_m$
whose joint distribution (given segment boundaries and $m=l+1$) is
\beq\label{ssd}
  \mbox{single segment:}\quad
  P(y_{ij},\mu_m|t_{m-1,m},1) \;=\; P(\mu_m)\prod_{t=i+1}^j P(y_t|\mu_m)
\eeq
by the model assumptions \req{lh} and \req{pc}. The probabilities
for a general but fixed segmentation are independent, i.e.
\bqa\label{PPP}
  P(y_{ij},\mu_{lm}|t_{lm},m-l)
  &=& \prod_{p=l+1}^m\left[P(\mu_p)\prod_{t=t_{p-1}+1}^{t_p} P(y_t|\mu_p)\right]
\\
  &=& P(y_{ih},\mu_{lp}|t_{lp},p-l)P(y_{hj},\mu_{pm}|t_{pm},m-p) \qmbox{(any $p$)}
\eqa
This is our key recursion. Consider now
\bqa\label{Qrec}
  Q(y_{ij},\mu_{lm}|m-l)
  &:=& ({\textstyle{j-i-1\atop m-l-1}})P(y_{ij},\mu_{lm}|t_l,t_m,m-l)
\\ \label{Qreca}
  &\stackrel{(a)}=& ({\textstyle{j-i-1\atop m-l-1}})
  \sum_{t_{lm}\,:\,i=t_l<...<t_m=j\hspace{-9ex}} P(y_{ij},\mu_{lm}|t_{lm},m-l)P(t_{lm}|m-l)
\\ \nonumber
  &\stackrel{(b)}=& \sum_{t_{lm}\,:\,i=t_l<...<t_m=j\hspace{-9ex}} P(y_{ij},\mu_{lm}|t_{lm},m-l)
\\ \nonumber
  &\stackrel{(c)}=& \sum_{t_p=i+p-l}^{j+p-m}
  \sum_{t_{lp}\,:\,i=t_l<...<t_p=h\hspace{-11ex}} P(y_{ih},\mu_{lp}|t_{lp},p-l)
  \sum_{t_{pm}\,:\,h=t_p<...<t_m=j\hspace{-11ex}} P(y_{hj},\mu_{pm}|t_{pm},m-p)
\\ \label{Qrecl}
  &=& \sum_{h=i+p-l}^{j+p-m} Q(y_{ih},\mu_{lp}|p-l)Q(y_{hj},\mu_{pm}|m-p)
\eqa
$(a)$ is just an instance of formula $P(A)=\sum_i P(A|H_i)P(H_i)$ for a
partitioning $(H_i)$ of the sample space. In $(b)$ we exploited
uniformity \req{ubp} of $P(t_{lm}|m-l)=({j-i-1\atop m-l-1})^{-1}$
and hence its independence from the concrete segmentation
$t_{lm}$. In $(c)$ we fix segment boundary $t_p$, sum over the
left and right segmentations, and finally over $t_p$.

\paradot{Left and right recursions}
If we integrate \req{Qrec} over $\mu_{lm}$, the integral
factorizes and we get a recursion in (a quantity that is
proportional to) the evidence of $y_{ij}$. Let us define more
generally $r^{th}$ ``Q-moments'' of $\mu'_t$.
\bqa\label{Qrrec}
  Q_t^r(y_{ij}|m-l) &:=& \int Q(y_{ij},\mu_{lm}|m-l) \mu'_t\!\,^r d\mu_{lm}
\\ \nonumber
  &=& \nq\sum_{h=i+p-l}^{t-1}\!\! Q^0(y_{ih}|p-l)Q_t^r(y_{hj}|p-l)
  +\!\!\! \sum_{h=t}^{j+p-m}\!\! Q_t^r(y_{ih}|m-p)Q^0(y_{hj}|m-p)
\eqa
Depending on whether $h<t$ or $h\geq t$, the $\mu'_t\!\,^r$ term
combines with the right or left $Q$ in recursion \req{Qrecl} to
$Q_t^r$, while the other $Q$ simply gets integrated to $Q_t^0=Q^0$
independent $t$. The recursion terminates with
\beq\label{Adef}
  A_{ij}^r \;:=\; Q_t^r(y_{ij}|1) \;=\;
  \int P(\mu_m)\prod_{t=i+1}^j P(y_t|\mu_m) \mu_m^r d\mu_m,
  \quad (0\leq i<j\leq n)
\eeq
Note $A_{ij}^0=P(y_{ij}|t_{m-1,m})$ is the evidence and
$A_{ij}^r/A_{ij}^0=\E[\mu_m^r|y_{ij},t_{m-1,m}]$ the $r^{th}$
moment of $\mu'_t=\mu_m$ in case $y_{ij}$ is modeled by a single
segment.
It is convenient to formally start the recursion with
$Q^0(y_{ij}|0)=\delta_{ij}=\{ {1 \text{ if } i=j\atop 0 \text{ else}\;\;\;}$ (consistent
with the recursion) with interpretation that (only) an empty data
set ($i=j$) can have 0 segments.
Since $p$ was an arbitrary split number, we can choose it conveniently.
We need a left recursion for $r=0$, $i=0$, $p-l=k$, and $m-p=1$:
\beqn
  L_{k+1,j} \;:=\; Q^0(y_{0j}|k+1)
  \;=\; \sum_{h=k}^{j-1} Q^0(y_{0h}|k)Q^0(y_{hj}|1)
  \;=\; \sum_{h=k}^{j-1} L_{kh}A^0_{hj}
\eeqn
That is (apart from binomial factors) the evidence of $y_{0j}$
with $k+1$ segments equals the evidence of $y_{0h}$ with $k$
segments times the single-segment evidence of $y_{hj}$, summed
over all locations $h$ of boundary $k$.
The recursion starts with $L_{1j}=A_{0j}^0$, or more conveniently
with $L_{0j}=\delta_{j0}$.
We also need a right recursion for
$r=0$, $j=n$, $p-l=1$, $m-p=k$:
\beqn
  R_{k+1,i} \;:=\; Q^0(y_{in}|k+1)
  \;=\; \sum_{h=i+1}^{n-k} Q^0(y_{ih}|1)Q^0(y_{hn}|k)
  \;=\; \sum_{h=i+1}^{n-k} A^0_{ih}R_{kh}
\eeqn
The recursion starts with $R_{1n}=A_{in}^0$, or more conveniently
with $R_{0i}=\delta_{in}$.

\paradot{Quantities of interest}
Note that
\beqn
  L_{kn} \;=\; R_{k0} \;=\; Q^0(\v y|k) \;=\; (\textstyle{n-1\atop k-1})P(\v y|k)
\eeqn
are proportional to the data evidence for fixed $k$. So the data
evidence can be computed as
\beq\label{E}
  E \;:=\; P(\v y) \;=\; \sum_{k=1}^n P(\v y|k)P(k)
  \;=\; {1\over k_{max}}\sum_{k=1}^{k_{max}} {L_{kn}\over({n-1\atop k-1})}
\eeq
The posterior of $k$ and its MAP estimate are
\beq\label{Ck}
  C_k \;:=\; P(k|\v y) \;=\; {P(\v y|k)P(k)\over P(\v y)} \;=\;
  {L_{kn}\over ({n-1\atop k-1})k_{max}E}
  \qmbox{and} \hat k=\mathop{\arg\max}\limits_{k=1..k_{max}} C_k
\eeq

\paradot{Segment boundaries}
We now determine the segment boundaries.
Consider recursion \req{Qrec} for $i=l=0$, $m=k$, $j=n$,
but keep $t_p=h$ fixed, i.e.\
do not sum over it.
Then \req{Qreca} and \req{Qrecl} reduce to the l.h.s.\ and r.h.s.\ of
\beq\label{PQQ}
  ({\textstyle{n-1\atop k-1}})P(\v y,\v\mu,t_p|k) \;=\;
  Q(y_{0h},\mu_{0p}|p)Q(y_{hn},\mu_{pk}|k-p)
\eeq
Integration over $\v\mu$ gives
\beqn
  ({\textstyle{n-1\atop k-1}})P(\v y,t_p|k) \;=\;
  Q^0(y_{0h}|p)Q^0(y_{hn}|k-p)
\eeqn
Hence the posterior probability that boundary $p$ is located at
$t_p=h$, given $\hat k$, is
\beq\label{Bph}
  B_{ph} \;:=\; P(t_p=h|\v y,\hat k) \;=\;
  { ({n-1\atop\hat k-1})P(\v y,t_p|\hat k)\over({n-1\atop\hat k-1})P(\v y|\hat k)}
  \;=\; {L_{ph}R_{\hat k-p,h}\over L_{\hat k n}}
\eeq
So our estimate for segment boundary $p$ is
\beq\label{hattp}
  \hat t_p \;:=\; \arg\max_h P(t_p=h|\v y,\hat k)
  \;=\; \arg\max_h\{ B_{ph} \}
  \;=\; \arg\max_h\{ L_{ph}R_{\hat k-p,h} \}
\eeq

\paradot{Segment levels}
Finally we need the segment levels, given the segment number $\hat k$ and
boundaries $\hat t$. The $r^{th}$ moment of segment $m$ with boundaries
$i=\hat t_{m-1}$ and $j=\hat t_m$ is
\beq\label{mumr}
  \widehat{\mu_m^r} \;=\; \E[\mu_m^r|\v y,\hat t,\hat k]
  \;=\; \E[\mu_m^r|y_{ij},\hat t_{m-1,m},1]
  \;=\; { \int P(y_{ij},\mu_m|\hat t_{m-1,m},1)\mu_m^r d\mu_m \over
          \int P(y_{ij},\mu_m|\hat t_{m-1,m},1) d\mu_m }
  \;=\; {A_{ij}^r\over A_{ij}^0}
\eeq
Note that this expression is independent of other segment boundaries and
their number, as it should.

\paradot{Regression curve}
Recursion \req{Qrrec} allows in principle to compute the
regression curve $\E[\mu'_t|\v y]$ by defining $(L_t^{r=1})_{kj}$
and $(R_t^{r=1})_{ki}$ analogous to $L_{kj}$ and $R_{ki}$, but
this procedure needs $O(n^3)$ space and $O(k_{max}n^3)$ time, one
$O(n)$ worse than our target performance.
We reduce probabilities of $\mu'_t$ to probabilities of $\mu_m$:
We exploit the fact that in every segmentation, $\mu'_t$ lies in
some segment. Let this (unique) segment be $m$ with (unique)
boundaries $i=t_{m-1}<t\leq t_m=j$. Then $\mu'_t=\mu_m$.
Summing now over all such segments we get
\beq\label{reg1}
  P(\mu'_t|\v y,k) \;=\;
  \sum_{m=1}^k\sum_{i=0}^{t-1}\sum_{j=t}^n P(\mu_m,t_{m-1}=i,t_m=j|\v y,k)
\eeq
By fixing $t_p$ in \req{Qreca} we arrived at \req{PQQ}.
Similarly, dividing the data into three parts and fixing $t_l$ and $t_m$
we can derive
\beqn
  ({\textstyle{n-1\atop k-1}})P(\v y,\v\mu,t_l,t_m|k)
  \;=\; Q(y_{0i},\mu_{0l}|l)Q(y_{ij}\mu_m|m-l)Q(y_{jn}\mu_{mk}|k-m)
\eeqn
Setting $l=m-1$, integrating over $\mu_{0l}$ and $\mu_{mk}$, dividing by
$({\textstyle{n-1\atop k-1}})P(\v y|k)$, and inserting into
\req{reg1}, we get
\beqn
  P(\mu'_t|\v y,k) \;=\;  {1\over L_{kn}}\sum_{m=1}^k\sum_{i<t\leq j}
  L_{m-1,i}Q(y_{ij},\mu_m|1)R_{k-m,j}
\eeqn
The posterior moments of $\mu'_t$, given $\hat k$, can hence be computed
by
\beq\label{muptr}
  \widehat{\mu'_t\!\,^r}
  \;=\; \sum_{i<t\leq j} F_{ij}^r
  \qmbox{with} F_{ij}^r:={1\over L_{\hat k n}}\sum_{m=1}^{\hat k} L_{m-1,i}A_{ij}^r R_{\hat k-m,j}
\eeq
While segment boundaries and values make sense only for fixed $k$
(we chose $\hat k$), the regression curve $\hat\mu'_t$ could
actually be averaged over all $k$ instead of fixing $k=\hat k$.

\paradot{Relative log-likelihood}
Another quantity of interest is how likely it is that $\v y$ is
sampled from $\hat f$. The log-likelihood of
$\v y$ is
\beqn
  ll \;:=\; \log P(\v y|\hat f)
      \;=\; \log P(\v y|\v{\hat\mu},\v{\hat t},\hat k)
      \;=\; \sum_{i=1}^n\log P(y_i|\hat\mu'_i,\s)
\eeqn
Like for the evidence, the number itself is hard to interpret. We
need to know how many standard deviations it is away from its
mean(=entropy). Since noise \req{lh} is i.i.d., mean and variance
of $ll$ are just $n$ times the mean and variance of the log-noise
distribution of a single data item. For Gaussian and Cauchy noise
we get
\bqan
  \mbox{Gauss:} & & \E[ll|\hat f] = \textstyle {n\over 2}\log(2\pi\e\hat\s^2),
         \qquad \Var[ll|\hat f] = {n\over 2}
\\
  \mbox{Cauchy:} & & \E[ll|\hat f] = \textstyle n\log(4\pi\hat\s), \hspace{2ex}
          \qquad \Var[ll|\hat f] = {n\over 3}\pi^2
\eqan

\section{Computing the Single Segment Distribution}\label{secSSD}

We now determine (at least in the Gaussian case efficient)
expressions for the moments \req{Adef} of the distribution
\req{ssd} of a single segment.

\paradot{Gaussian model}
For Gaussian noise \req{Gn} and prior \req{Gp} we get
\beqn
  A_{ij}^r \;=\;
  \left({1\over\sqrt{2\pi}\s}\right)^d {1\over\sqrt{2\pi}\rho}
  \int_{-\infty}^\infty \mbox{\large e}^{\textstyle-{1\over 2\s^2}\sum_{t=i+1}^j(y_t\!-\!\mu_m)^2
   -{1\over 2\rho^2}(\mu_m\!-\!\nu)^2} \mu_m^r\, d\mu_m
\eeqn
where $d=j-i$. This is an unnormalized Gaussian integral with the
following normalization, mean, and variance \cite[Sec.10.2]{Bolstad:04}:
\bqa\label{GyP}
  P(y_{ij}|t_{m-1,m}) &=&
  A_{ij}^0
  \;=\;
   {\exp\!\big\{ {1\over 2\s^2}\big[{(\sum_t(y_t\!-\!\nu))^2\over d+\s^2/\rho^2} \!-\! \sum_t(y_t\!-\!\nu)^2\big]\big\}
   \over(2\pi\s^2)^{d/2}(1\!+\!d\rho^2/\s^2)^{1/2}}
\\ \label{Gmm}
  \E[\mu_m|y_{ij},t_{m-1,m}]
  &=& {A_{ij}^1\over A_{ij}^0}
  \;=\; {\rho^2(\sum_t y_t)+\s^2\nu \over d\rho^2+\s^2}
  \;\approx\; {1\over d}\sum_t y_t
\\ \label{Gmvar}
  \Var[\mu_m|y_{ij},t_{m-1,m}]
  &=& {A_{ij}^2\over A_{ij}^0} - \Big({A_{ij}^1\over A_{ij}^0}\Big)^2
  \;=\; \Big[{d\over\s^2}+{1\over\rho^2}\Big]^{-1}
  \;\approx\; {\s^2\over d}
\eqa
where $\Sigma_t$ runs from $i+1$ to $j$.
The mean/variance is just the weighted average of the
mean/variance of $y_{ij}$ and $\mu_m$.
One may prefer to use the segment prior only for determining
$A_{ij}^0$, but use the unbiased estimators ($\approx$) for the moments.
Higher moments $A_{ij}^r$ can also be computed from the central
moments
\beqn
  \E[(\mu_m-A_{ij}^1/A_{ij}^0)^r|y_{ij},t_{m-1,m}]
  \;=\; {1\!\cdot\!3\!\cdot...\cdot\!(r-1)\over[d\s^{-2}+\rho^{-2}]^{r/2}}
  \;\approx\; 1\!\cdot\!3\!\cdot...\cdot\!(r-1)\!\cdot\!
  \Big({\s^2\over d}\Big)^{r/2}
\eeqn
for even $r$, and 0 for odd $r$.

\paradot{Other models}
Analytic expressions for $A_{ij}^r$ are possible for all
distributions in the exponential family. For others like Cauchy
we need to perform integral \req{Adef} numerically.
A very simple approximation is to replace the integral by a sum on
a uniform grid: The stepsize/range of the grid should be some
fraction/multiple of the typical scale of the integrand, and the
center of the grid should be around the mean. A crude estimate of
the mean and scale can be obtained from the Gaussian model
\req{Gmm} and \req{Gmvar}.
Or even simpler, use the estimated global mean and variance
\req{eGmv}, and in-segment variance \req{sEst} for determining the
range (e.g.\ $[\hat\nu-25\hat\rho,...,\hat\nu+25\hat\rho]$) and
stepsize (e.g.\ $\hat\s/10$) of one grid used for all $A_{ij}^r$.
Note that if $y_{ij}$ really stem from one segment, the integrand
is typically unimodal and the above estimates for stepsize and
range are reasonable, hence the approximation will be good. If
$y_{ij}$ ranges over different segments, the discretization may be
crude, but since in this case, $A_{ij}^r$ is (very) small, crude
estimates are sufficient.
Note also that even for the heavy-tailed Cauchy distribution, the
first and second moments $A_{ij}^1$ and $A_{ij}^2$ exist, since
the integrand is a product of at least two Cauchy distributions,
one prior and one noise for each $y_t$.
Preferably, standard numerical integration routines (which are
faster, more robust and more accurate) should be used.

\section{Determination of the Hyper-Parameters}\label{secHP}

\paradot{Hyper-Bayes and Hyper-ML}
The developed regression model still contains three
(hyper)parameters, the global variance $\rho^2$ and mean $\nu$ of
$\v\mu$, and the in-segment variance $\s^2$. If they are not known, a
proper Bayesian treatment would be to assume a hyper-prior over
them and integrate them out.
Since we do not expect a significant influence of the hyper-prior
(as long as chosen reasonable) on the quantities of interest, one
could more easy proceed in an empirical Bayesian way and choose
the parameters such that the evidence $P(\v y|\s,\nu,\rho)$ is
maximized (``hyper-ML''). (We restored the till now omitted
dependency on the hyper-parameters).

Exhaustive (grid) search for the hyper-ML parameters is expensive.
For data which is indeed noisy piecewise constant, $P(\v
y|\s,\nu,\rho)$ is typically unimodal\footnote{A little care is
necessary with the in-segment variance $\s^2$. If we set it (extremely
close) to zero, all segments will consist of a single data point
$y_i$ with (close to) infinite evidence (see e.g.\ \req{GyP}).
Assuming $k_{max}<n$ eliminates this unwished maximum. Greedy
hill-climbing with proper initialization will also not be
fooled.} in $(\s,\nu,\rho)$ and the global maximum can be found more
efficiently by greed hill-climbing, but even this may cost a
factor of 10 to 1000 in efficiency. Below we present a very simple
and excellent heuristic for choosing $(\s,\nu,\rho)$.

\paradot{Estimate of global mean and variance $\nu$ and $\rho$}
A reasonable choice for the level mean and variance $\nu$ and $\rho$
are the empirical global mean and variance of the data $\v y$.
\beq\label{eGmv}
  \hat\nu \;\approx\; {1\over n} \sum_{t=1}^n y_t \qmbox{and}
  \hat\rho^2 \;\approx\; {1\over n-1} \sum_{t=1}^n (y_t-\hat\nu)^2
\eeq
This overestimates the variance $\rho^2$ of the segment levels,
since the expression also includes the in-segment variance $\s^2$,
which one may want to subtract from this expression.

\paradot{Estimate of in-segment variance $\s^2$}
At first there seems little hope of estimating the in-segment
variance $\s^2$ from $\v y$ without knowing the segmentation, but
actually we can use a simple trick. If $\v y$ would belong to a
single segment, i.e.\ the $y_t$ were i.i.d.\ with variance $\s^2$,
then the following expressions for $\s^2$ would hold:
\beqn
  \E[{1\over n}\sum_{t=1}^n(y_t-\mu_1)^2]
  \;=\; \s^2
  \;=\; {1\over 2(n-1)}\E[\sum_{t=1}^{n-1}(y_{t+1}-y_t)^2]
\eeqn
i.e.\ instead of estimating $\s^2$ by the squared deviation of the
$y_t$ from their mean, we can also estimate $\s^2$ from the
average squared difference of successive $y_t$. This remains true
even for multiple segments if we exclude the segment boundaries
in the sum. On the other hand, if the number of segment boundaries is
small, the error from including the boundaries will be small, i.e.\ the
second expression remains approximately valid.
More precisely, we have within a segment and at the boundaries
\beqn
  \E\nq\;\sum_{t=t_{m-1}+1}^{t_m-1}\nq\;(y_{t+1}-y_t)^2 = 2(t_m-t_{m-1}-1)\s^2
  \qmbox{and}
  \E(y_{t_m+1}-y_{t_m})^2 = 2\s^2+(\mu_{m+1}-\mu_m)^2
\eeqn
Summing over all $k$ segments and boundaries and solving w.r.t.\
$\s^2$ we get
\bqan
  \s^2 &=& {1\over 2(n-1)}\left\{\E\bigg[\sum_{t=1}^{n-1}(y_{t+1}-y_t)^2\bigg]
  - \sum_{m=1}^{k-1}(\mu_{m+1}-\mu_m)^2\right\}
\\
  &=& {1\over 2(n-1)}\E\bigg[\sum_{t=1}^{n-1}(y_{t+1}-y_t)^2\bigg]\cdot
  \bigg[1-O\Big({k\over n}{\rho^2\over\s^2}\Big)\bigg]
\eqan
The last expression holds, since there are $k$ boundaries in $n$
data items, and the ratio between the variance of $\v\mu$ to the
in-segment variance is $\rho^2/\s^2$. Hence we may estimate
$\s^2$ by the upper bound
\beq\label{sEst}
  \hat\s^2 \;\approx\; {1\over 2(n-1)}\sum_{t=1}^{n-1}(y_{t+1}-y_t)^2
\eeq
If there are not too many segments ($k\ll n$) and the regression
problem is hard (high noise $\rho\lesssim\s$), this is a very good
estimate. In case of low noise ($\rho\gg\s$), regression is very
easy, and a crude estimate of $\s^2$ is sufficient. If there are
many segments, $\hat\s^2$ tends to overestimate $\s^2$, resulting
in a (marginal) bias towards estimating fewer segments (which is
then often welcome).

If the estimate is really not sufficient, one may use \req{sEst}
as an initial estimate for determining an initial segmentation
$\hat t$, which then can be used to compute an improved estimate
of $\hat\s^2$, and possibly iterate.

\paradot{Hyper-ML estimates}
Expressions \req{eGmv} are the standard estimates of mean and
variance of a distribution. They are particularly suitable for
(close to) Gaussian distributions, but also for others, as long as
$\nu$ and $\rho$ parameterize mean and variance. If
mean and variance do not exist or the distribution is quite
heavy-tailed, we need other estimates.
The ``ideal'' hyper-ML estimates may be approximated as follows.
If we assume that each data point lies in its own segment,
we get
\beqn
  (\hat\nu,\hat\rho) \;\approx\;
  \mathop{\arg\max}\limits_{(\nu,\rho)}
  \prod_{t=1}^n P(y_t|\hat\s,\nu,\rho) \qmbox{with}
\eeqn
\beq\label{Pysnr}
  P(y_t|\s,\nu,\rho) \;=\; \int P(y_t|\mu,\s)P(\mu|\nu,\rho)d\mu
\eeq
The in-segment variance $\hat\s^2$ can be estimated similarly to the
last paragraph considering data differences and ignoring segment
boundaries:
\beqn
  \hat\s \;\approx\; \arg\max_\s \prod_{t=1}^{n-1} P(y_{t+1}-y_t|\s)
  \qmbox{with}
\eeqn
\beq\label{PDys}
  P(y_{t+1}-y_t=\Delta|\s) \;\approx\; \int_{-\infty}^\infty
  P(y_{t+1}=a+\Delta|\mu,\s) P(y_t=a|\mu,\s)da
\eeq
Note that the last expression is independent of the segment level
(this was the whole reason for considering data differences) and
exact iff $y_t$ and $y_{t+1}$ belong to the same segment.
In general (beyond the exponential family)
$(\hat\nu,\hat\rho,\hat\s)$ can only be determined numerically.

\paradot{Using median and quartile}
We present some simpler estimates based on median and quartiles.
Let $[\v y]$ be the data vector $\v y$, but sorted in ascending
order. Then, item $[\v y]_{\alpha n}$ (where the index is assumed
to be rounded up to the next integer) is the $\alpha$-quantile of
empirical distribution $\v y$. In particular $[\v y]_{n/2}$ is the
median of $\v y$. It is a consistent (and robust to outliers)
estimator of the mean segment level
\beq\label{nuEstm}
  \hat\nu \;\approx\; [\v y]_{n/2}
\eeq
if noise and segment levels have symmetric distributions. Further,
half of the data points lie in the interval $[a,b]$, where $a:=[\v
y]_{n/4}$ is the first and $b:=[\v y]_{3n/4}$ is the last quartile
of $\v y$. So, using \req{Pysnr}, $\hat\rho$ should be estimated
such that
\beqn
  P(a\leq y_t\leq b|\s,\hat\nu,\hat\rho) \;\stackrel!\approx \odt
\eeqn
Ignoring data noise (assuming $\s\approx 0$), we get
\beq\label{rhoEstm}
  \hat\rho \;\approx\; {[\v y]_{3n/4}-[\v y]_{n/4}\over 2\alpha}
  \qmbox{with $\alpha=1$ for Cauchy and $\alpha\doteq 0.6744$ for Gauss,}
\eeq
where $\alpha$ is the quartile of the standard
Cauchy/Gauss/other segment prior. For the data noise $\s$ we
again consider the differences $\Delta_t:=y_{t+1}-y_t$. Using
\req{PDys}, $\hat\s$ should be estimated such that
\beqn
  P(a'\leq y_{t+1}-y_t\leq b'|\hat\s) \stackrel!\approx \odt
\eeqn
where $a'=[\v\Delta]_{n/4}$ and $b'=[\v\Delta]_{3n/4}\approx -a'$.
One can show that
\beq\label{sEstm}
  \hat\s \;\approx\; {[\v\Delta]_{3n/4}-[\v\Delta]_{n/4}\over 2\beta}
  \qmbox{with $\beta=2$ for Cauchy and $\beta\doteq 0.6744\sqrt{2}$ for Gauss,}
\eeq
where $\beta$ is the quartile of the one time with itself
convolved standard Cauchy/Gauss/other (noise) distribution.
Use of quartiles for estimating $\s$ is robust to the ``outliers''
caused by the segment boundaries, so yields better estimates than
\req{sEst} if noise is low.
Again, if the estimates are really not sufficient, one may
iteratively improve them.

\section{The Algorithm}\label{secAlg}

The computation of $A$, $L$, $R$, $E$, $C$, $B$, $\hat t_p$,
$\widehat{\mu_m^r}$, $F$, and $\widehat{\mu'_t\!\,^r}$ by the
formulas/recursions derived in Section \ref{secES}, are
straightforward. In \req{Adef} one should compute the product, or
in \req{GyP}, \req{Gmm}, \req{Gmvar} the sum, incrementally from
$j\leadsto j+1$. Similarly $\widehat{\mu'_t\!\,^r}$ should be
computed incrementally by
\beqn
  \widehat{\mu'_{t+1}\nq\,^r} \;=\;
  \widehat{\mu'_t\!\,^r} - \sum_{i=0}^{t-1}F_{it}^r + \sum_{j=t+1}^n F_{tj}^r
\eeqn
Typically $r=0,1,2$. In this way, all quantities can be computed
in time $O(k_{max}n^2)$ and space $O(n^2)$. Space can be reduced
to $O(k_{max}n)$ by computing $A$ on-the-fly in the various
expressions at the cost of a slowdown by a constant factor.
Table \ref{PCRegAlg} contains the algorithm in pseudo-C code. The
complete code including examples and data is available at
\cite{Hutter:05pcrcode}.
Since $A^0$, $L$, $R$, and $E$ can be exponentially large in $n$,
i.e.\ huge or tiny, actually their logarithm has to be computed
and stored. In the expressions, the logarithm is pulled in by
$\log(x\cdot y)=\log(x)+\log(y)$ and
$\log(x+y)=\log(x)+\log(1+\exp(\log(y)-\log(x))$ for $x>y$ and
similarly for $x<y$. Instead of $A_{ij}^r$ we have to compute
$A_{ij}^r/A_{ij}^0$ by pulling the denominator into the integral.

\begin{table} \def\algitsep{\itemsep=0ex}
\vspace{-3.5ex}
{\bf\caption{\label{PCRegAlg}Regression algorithm in pseudo C code\rm}}
\vspace{1ex}
{\bf\boldmath EstGauss($\v y,n$)} and {\bf\boldmath EstGeneral($\v
y,n,\alpha,\beta$)} compute from data $(y_1,...,y_n)$, estimates
for $\nu$, $\rho$, $\s$ (hat `$\hat{\rule{1ex}{0ex}}$' omitted),
and from that the evidence $A_{ij}^0$ of a single segment ranging
from $i+1$ to $j$, and corresponding first and second moments
$A_{ij}^1$ and $A_{ij}^2$. The expressions \req{eGmv}, \req{sEst},
\req{GyP}, \req{Gmm}, \req{Gmvar} are used in EstGauss() for
Gaussian noise and prior, and \req{nuEstm}, \req{rhoEstm},
\req{sEstm} and numerical integration on a uniform Grid in
EstGeneral() for arbitrary noise and prior $P$, e.g.\ Cauchy. $[\v
y]$ denotes the sorted $\v y$ array, Grid is the uniform
integration grid, $+=$ and $*=$ are additive/multiplicative
updates, and $\aai$ denotes arrays.

\vspace*{1ex}
\begin{minipage}[t]{0.53\textwidth}
\begin{list}{}{\parskip=0ex\parsep=0ex\algitsep\leftmargin=0ex\labelwidth=0ex}
  \item {\bf\boldmath EstGauss($\v y,n$)}
  \begin{list}{}{\parskip=0ex\parsep=0ex\algitsep\leftmargin=2ex\labelwidth=1ex\labelsep=1ex}
    \item[$\lceil$] $\nu={1\over n}\sum_{t=1}^n y_t$;
    \item $\rho^2={1\over n-1}\sum_{t=1}^n (y_t-\nu)^2$;
    \item $\s^2={1\over 2(n-1)}\sum_{t=1}^{n-1}(y_{t+1}-y_t)^2$;
    \item for($i=0..n$)
    \begin{list}{}{\parskip=0ex\parsep=0ex\algitsep\leftmargin=2ex\labelwidth=1ex\labelsep=1ex}
      \item[$\lceil$] $m=0$; $s=0$;
      \item for($j=i+1..n$)
      \begin{list}{}{\parskip=0ex\parsep=0ex\algitsep\leftmargin=2ex\labelwidth=1ex\labelsep=1ex}
        \item[$\lceil$] $d=j-i$; $m+=y_j-\nu$; $s+=(y_j-\nu)^2$;
        \item $A_{ij}^0=\displaystyle{\exp\{ {1\over 2\s^2}[{m^2\over d+\s^2/\rho^2}-s]\}
                        \over (2\pi\s^2)^{d/2}(1+d\rho^2/\s^2)^{1/2}}$;
        \item $A_{ij}^1=A_{ij}^0(\nu+m/d)$; 
      \end{list}
      \item[$\lfloor$] $\lfloor$ $A_{ij}^2 = A_{ij}^0 ((A_{ij}^1/A_{ij}^0)^2 + \s^2/d)$;
    \end{list}
    \item [$\lfloor$] {\bf\boldmath return ($A_{\aai\aai}^\aai,\nu,\rho,\s$); }
  \end{list}
\end{list}
\end{minipage}
\begin{minipage}[t]{0.46\textwidth}
\begin{list}{}{\parskip=0ex\parsep=0ex\algitsep\leftmargin=0ex\labelwidth=0ex}
  \item {\bf\boldmath EstGeneral($\v y,n,\alpha,\beta$)}
  \begin{list}{}{\parskip=0ex\parsep=0ex\algitsep\leftmargin=2ex\labelwidth=1ex\labelsep=1ex}
    \item[$\lceil$] $\nu=[\v y]_{n/2}$;
    \item $\rho=([\v y]_{3n/4}-[\v y]_{n/4})/2\alpha$;
    \item for($t=1..n-1$) $\Delta_t=y_{t+1}-y_t$;
    \item $\sigma=([\v\Delta]_{3n/4}-[\v\Delta]_{n/4})/2\beta$;
    \item Grid$=({\s\over 10}\SetZ)\cap[\nu-25\rho,\nu+25\rho]$;
    \item for($i=0..n$)
    \begin{list}{}{\parskip=0ex\parsep=0ex\algitsep\leftmargin=2ex\labelwidth=1ex\labelsep=1ex}
      \item[$\lceil$] for($\mu\in$Grid) $R_\mu=P(\mu|\nu,\rho)$;
      \item for($j=i+1..n$)
      \begin{list}{}{\parskip=0ex\parsep=0ex\algitsep\leftmargin=2ex\labelwidth=1ex\labelsep=1ex}
        \item[$\lceil$] for($\mu\in$Grid) $R_\mu*=P(y_j|\mu,\s)$;
      \end{list}
      \item[$\lfloor$] $\lfloor$ $A_{ij}^r={\s\over 10}\sum_{\mu\in\text{Grid}}R_\mu\,\mu^r$; $(r=0,1,2)$
    \end{list}
    \item [$\lfloor$] {\bf\boldmath return ($A_{\aai\aai}^\aai,\nu,\rho,\s$); }
  \end{list}
\end{list}
\end{minipage}

\vspace{-1ex}
\begin{minipage}{0.44\textwidth}
{\bf\boldmath Regression($\v A,n,k_{max}$)} takes $\v A$, $n$, and an upper bound
on the number of segments $k_{max}$, and computes the evidence $E=P(\v y)$ \req{E}, %
the probability $C_k=P(k|\v y)$ of $k$ segments %
and its MAP estimate $\hat k$ \req{Ck}, %
the probability $B_i=P(\exists p:t_p=i|\v y,\hat k)$ that a boundary is at $i$ \req{Bph} %
and the MAP location $\hat t_p$ of the $p^{th}$ boundary \req{hattp}, %
the first and second segment level moments $\mu_p$ and $\mu_p^2$ of all segments $p$ \req{mumr},
and the Bayesian regression curve $\mu'_t$ and its second moment $\mu'_t\!\,^2$ \req{muptr}.
\end{minipage}
\hspace{0.02\textwidth}
\begin{minipage}{0.54\textwidth}
\begin{list}{}{\parskip=0ex\parsep=0ex\algitsep\leftmargin=0ex\labelwidth=0ex}
  \item {\bf\boldmath Regression($\v A_{\aai\aai}^\aai,n,k_{max}$)}
  \begin{list}{}{\parskip=0ex\parsep=0ex\algitsep\leftmargin=2ex\labelwidth=1ex\labelsep=1ex}
    \item[$\lceil$] for($i=0..n$) \{ $L_{0i}=\delta_{i0}$; $R_{0i}=\delta_{in}$; \}
    \item for($k=0..n-1$)
    \begin{list}{}{\parskip=0ex\parsep=0ex\algitsep\leftmargin=2ex\labelwidth=1ex\labelsep=1ex}
      \item[$\lceil$] for($i=0..n$) $L_{k+1,i}=\sum_{h=k}^{i-1}L_{kh}A^0_{hi}$;
      \item[$\lfloor$] for($i=0..n$) $R_{k+1,i}=\sum_{h=i+1}^{n-k}A^0_{ih}R_{kh}$;
    \end{list}
    \item $E=k_{max}^{-1}\sum_{k=1}^{k_{max}} L_{kn}/({n-1\atop k-1})$;
    \item for($k=0..k_{max}$) $C_k= L_{kn}/[({n-1\atop k-1})k_{max}E]$;
    \item $\hat k=\mathop{\arg\max}_{k=1..k_{max}}\{C_k\}$;
    \item for($i=0..n$) $B_i=\sum_{p=0}^{\hat k}L_{pi}R_{\hat k-p,i}/L_{\hat k n}$;
    \item for($p=0..\hat k$) $\hat t_p=\arg\max_h\{ L_{ph}R_{\hat k-p,h} \}$;
    \item for($p=1..\hat k$) $\widehat{\mu_p^r}=A_{\hat t_{p-1}\hat t_p}^r/A_{\hat t_{p-1}\hat t_p}^0$; $(r=1,2)$
    \item for($i=0..n)$ for($j=i+1..n$)
    \item $[$ $F_{ij}^r=\sum_{m=1}^{\hat k} L_{m-1,i}A_{ij}^r R_{\hat k-m,j}/L_{\hat k n}$;
    \item $\mu'_0\!^r=0$; $(r=1,2)$
    \item for($t=0..n-1$)
    \item $[$ $\widehat{\mu'_{t+1}\nq\,^r} \;=
           \widehat{\mu'_t\!\,^r}-\sum_{i=0}^{t-1}F_{it}^r + \sum_{j=t+1}^n F_{tj}^r$
    \item [$\lfloor$] {\bf\boldmath return ($E,C_\aai,\hat k,B_\aai,\hat t_\aai,\widehat{\mu_\aai^r},\widehat{\mu_\aai'\!\,^r}$); }
  \end{list}
\end{list}
\end{minipage}

\end{table}

\section{Synthetic Examples}\label{secSE}

\begin{figure}
\begin{minipage}[t]{0.49\textwidth}
\includegraphics[width=\textwidth]{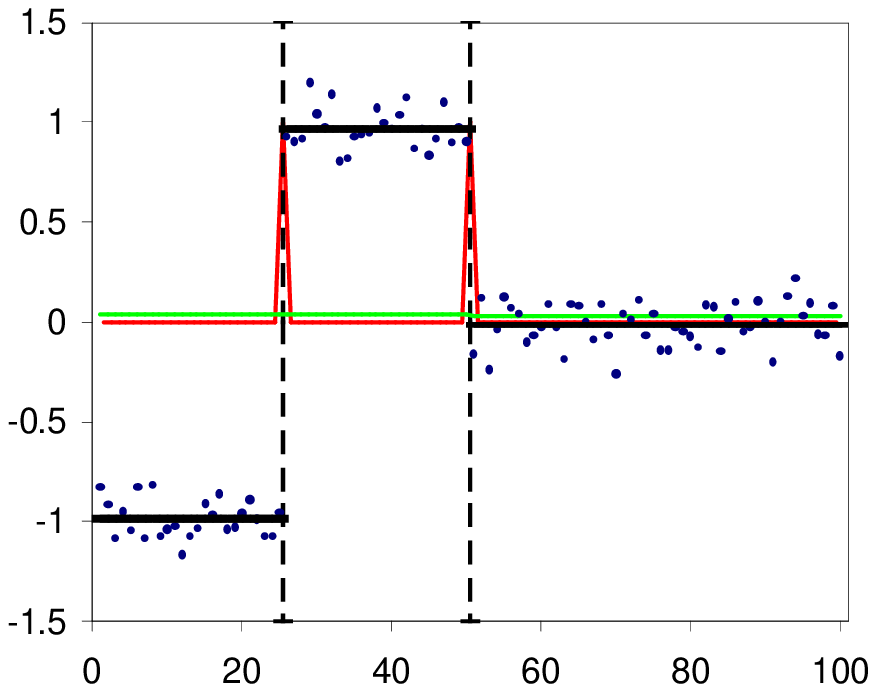}
\caption{\label{figGLPCR}[GL: low Gaussian noise]
data (blue), PCR (black), BP (red), and variance$^{1/2}$ (green).}
\end{minipage}
\hspace{0.02\textwidth}
\begin{minipage}[t]{0.49\textwidth}
\includegraphics[width=\textwidth]{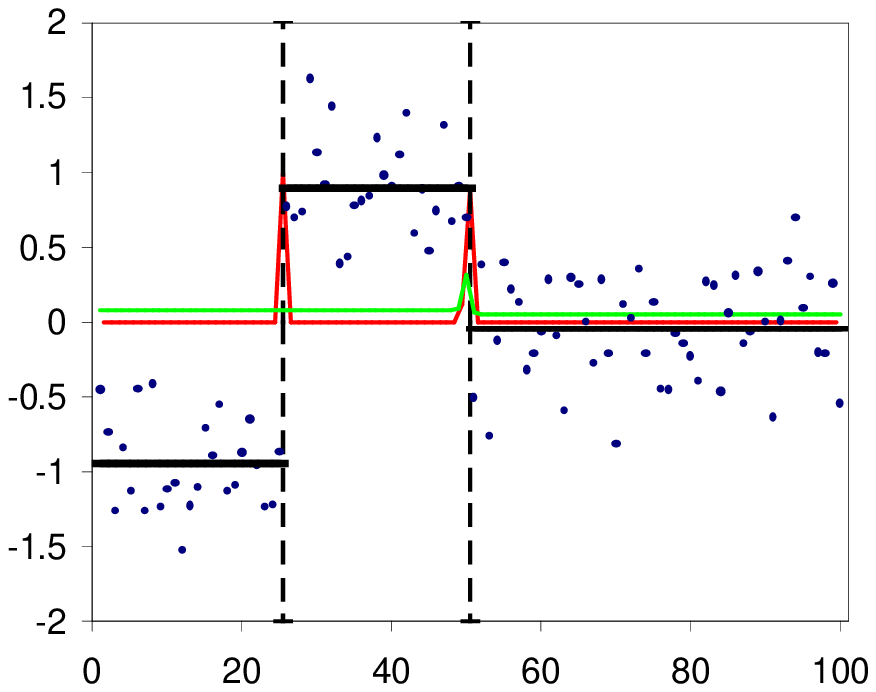}
\caption{\label{figGMPCR}[GM: medium Gaussian noise]
data (blue), PCR (black), BP (red), and variance$^{1/2}$ (green).}
\end{minipage}
\end{figure}

\begin{figure}
\begin{minipage}[t]{0.49\textwidth}
\includegraphics[width=\textwidth]{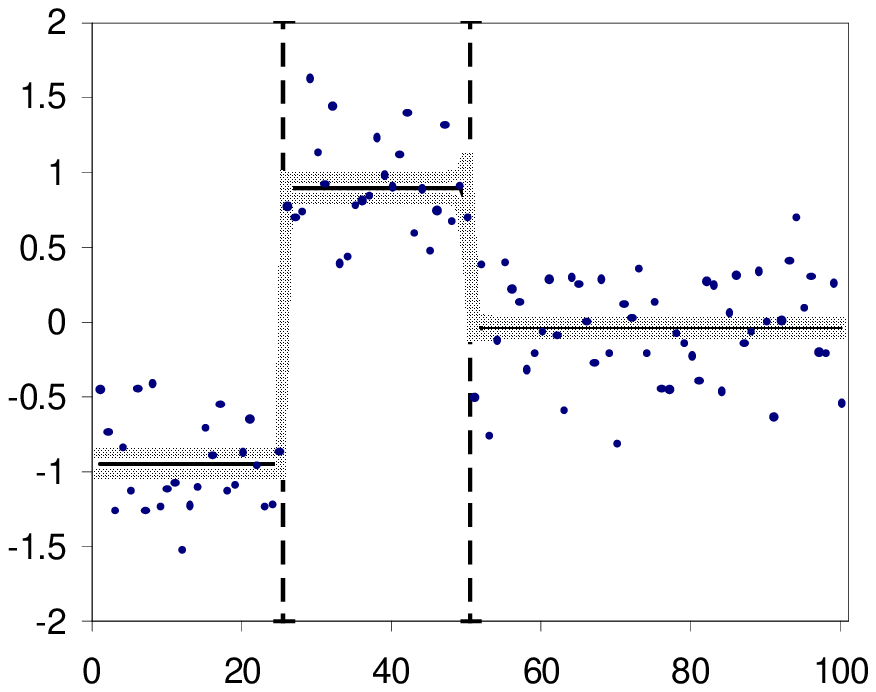}
\caption{\label{figGMBR}[GM: medium Gaussian noise]
data with Bayesian regression $\pm$ 1 std.-deviation.}
\end{minipage}
\hspace{0.02\textwidth}
\begin{minipage}[t]{0.49\textwidth}
\includegraphics[width=\textwidth]{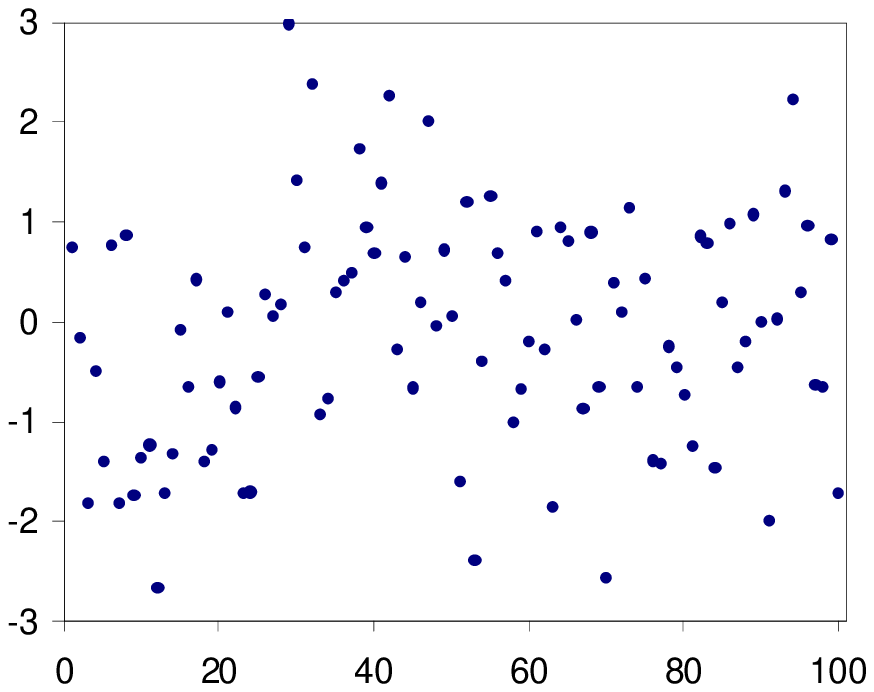}
\caption{\label{figGHData}[GH: high Gaussian noise]
data.}
\end{minipage}
\end{figure}

\begin{figure}
\begin{minipage}[t]{0.49\textwidth}
\includegraphics[width=\textwidth]{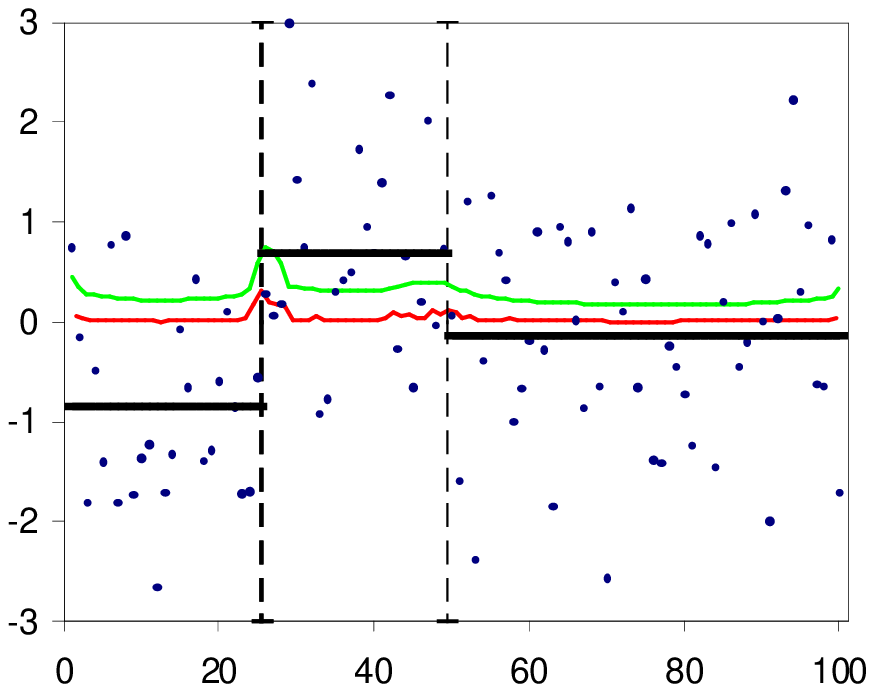}
\caption{\label{figGHPCR}[GH: high Gaussian noise]
data (blue), PCR (black), BP (red), and variance$^{1/2}$ (green).}
\end{minipage}
\hspace{0.02\textwidth}
\begin{minipage}[t]{0.49\textwidth}
\includegraphics[width=\textwidth]{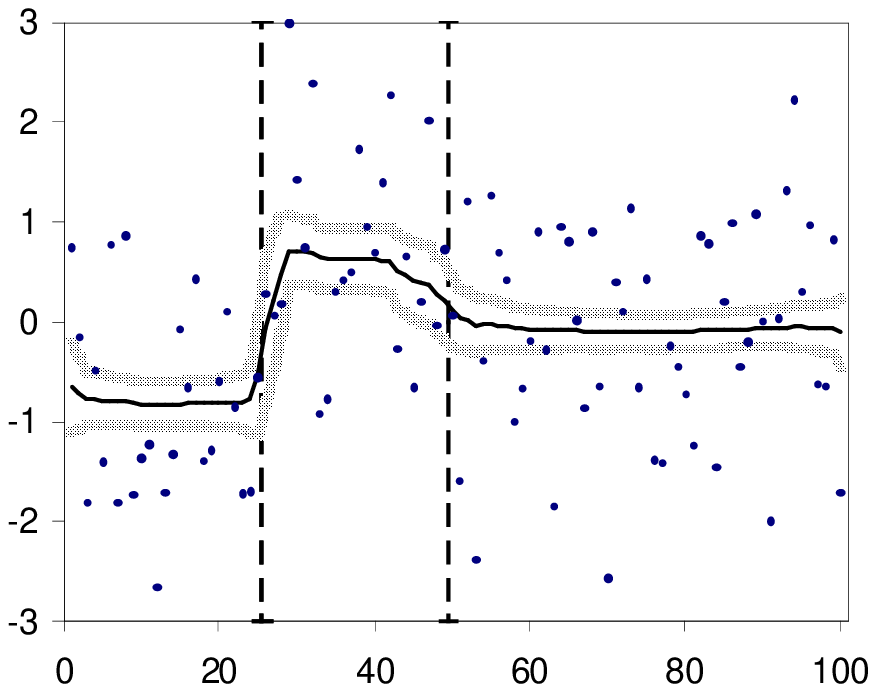}
\caption{\label{figGHBR}[GH: high Gaussian noise]
data with Bayesian regression $\pm$ 1 std.-deviation.}
\end{minipage}
\end{figure}

\begin{figure}
\begin{minipage}[t]{0.49\textwidth}
\includegraphics[width=\textwidth]{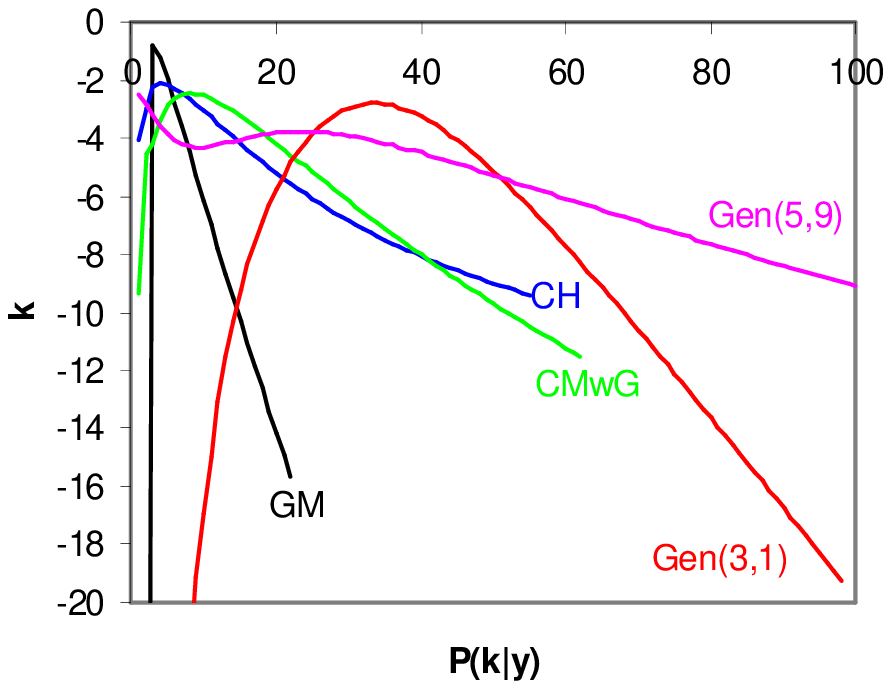}
\caption{\label{figPKY}Posterior segment number probability
$P(k|\v y)$ for medium Gaussian noise (GM, black), high Cauchy
noise (CH, blue), medium Cauchy noise with Gaussian regression
(CMwG, green), aberrant gene copy \# of chromosome 1 (Gen(3,1),
red), normal gene copy \# of chromosome 9 (Gen(5,9), pink).}
\end{minipage}
\hspace{0.02\textwidth}
\begin{minipage}[t]{0.49\textwidth}
\includegraphics[width=\textwidth]{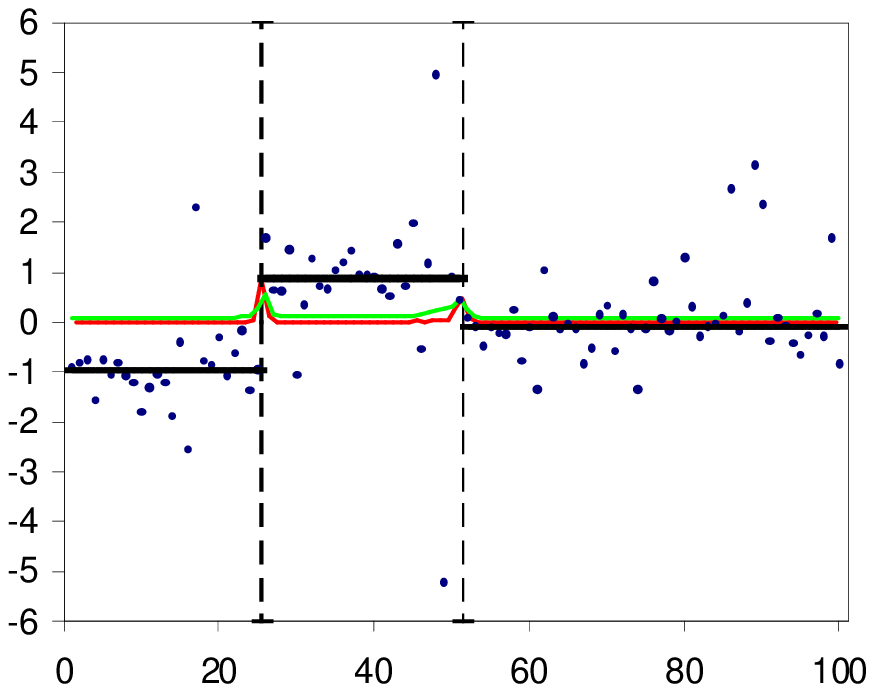}
\caption{\label{figCMPCR}[CM: medium Cauchy noise]
data (blue), PCR (black), BP (red), and variance$^{1/2}$ (green).}
\end{minipage}
\end{figure}

\begin{figure}
\begin{minipage}[t]{0.49\textwidth}
\includegraphics[width=\textwidth]{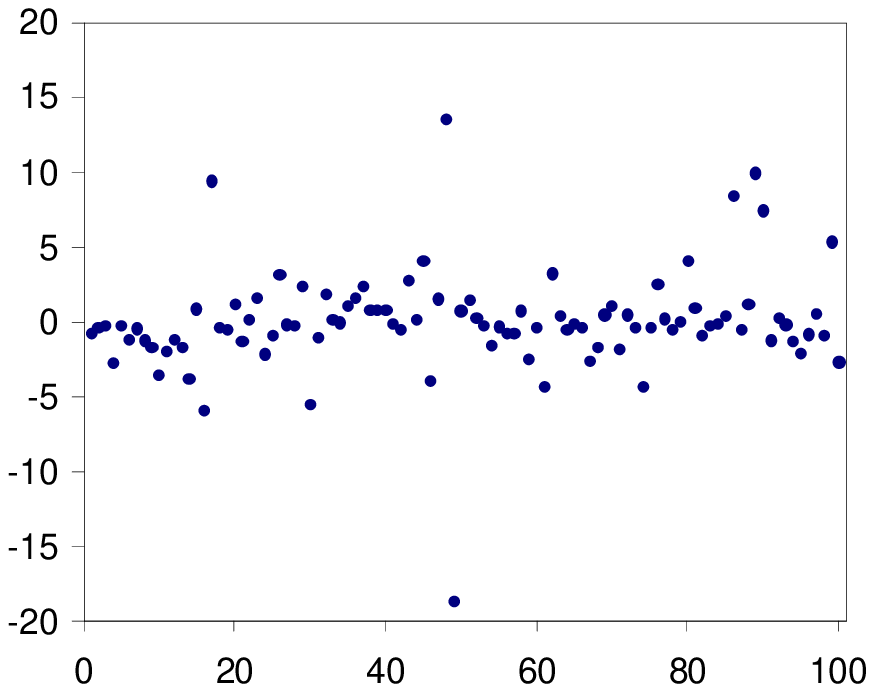}
\caption{\label{figCHData}[CH: high Cauchy noise]
data.}
\end{minipage}
\hspace{0.02\textwidth}
\begin{minipage}[t]{0.49\textwidth}
\includegraphics[width=\textwidth]{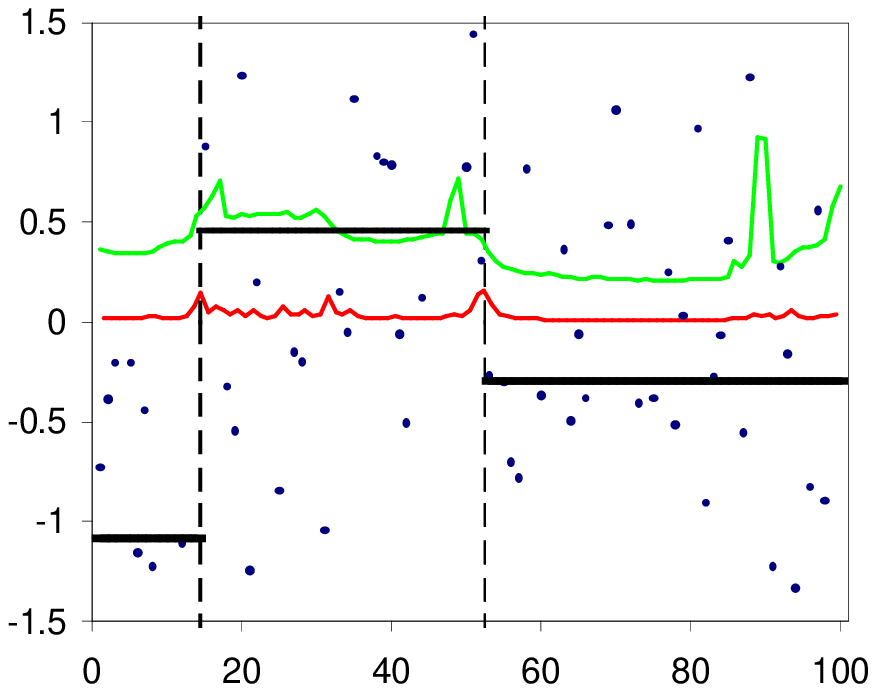}
\caption{\label{figCHPCR}[CH: high Cauchy noise]
data (blue), PCR (black), BP (red), and variance$^{1/2}$ (green).}
\end{minipage}
\end{figure}

\begin{figure}
\begin{minipage}[t]{0.49\textwidth}
\includegraphics[width=\textwidth]{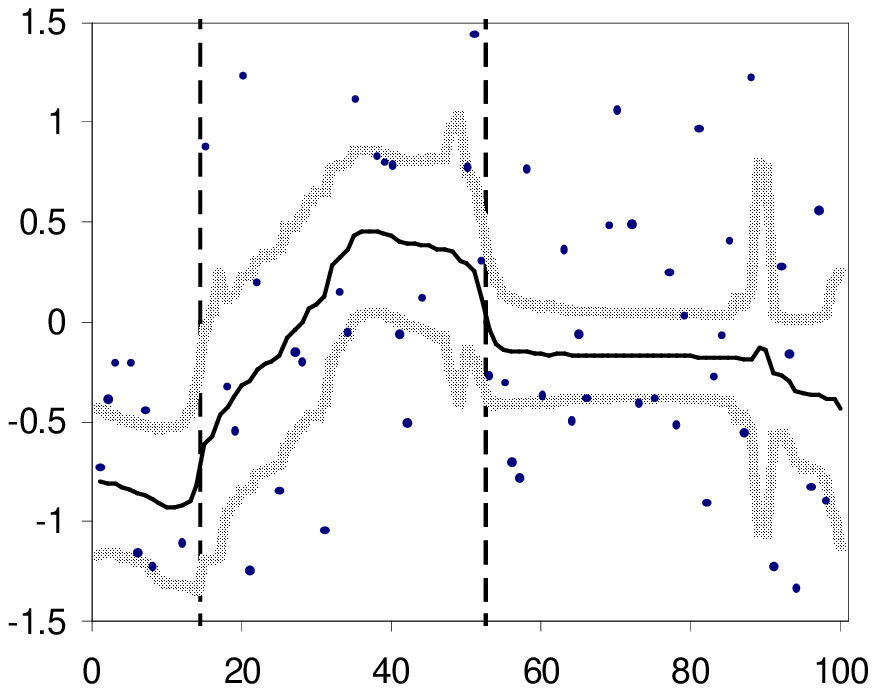}
\caption{\label{figCHBR}[CH: high Cauchy noise]
data with Bayesian regression $\pm$ 1 std.-deviation.}
\end{minipage}
\hspace{0.02\textwidth}
\begin{minipage}[t]{0.49\textwidth}
\includegraphics[width=\textwidth]{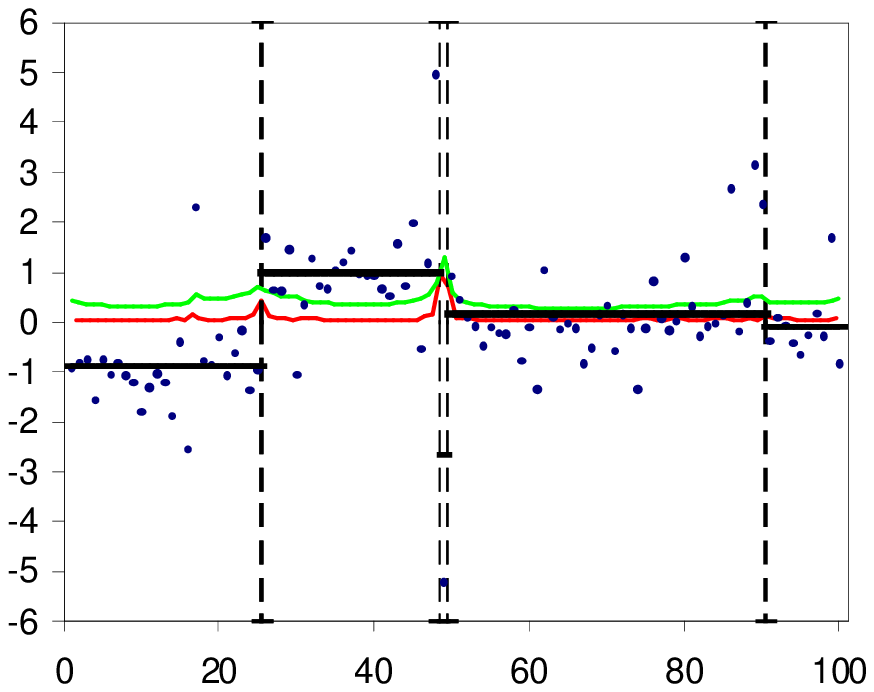}
\caption{\label{figCMPCRwG}[CMwG: medium Cauchy noise]
data (blue), but with Gaussian PCR (black), BP (red),
and variance$^{1/2}$ (green).}
\end{minipage}
\end{figure}

\begin{figure}
\begin{minipage}[t]{0.49\textwidth}
\includegraphics[width=\textwidth]{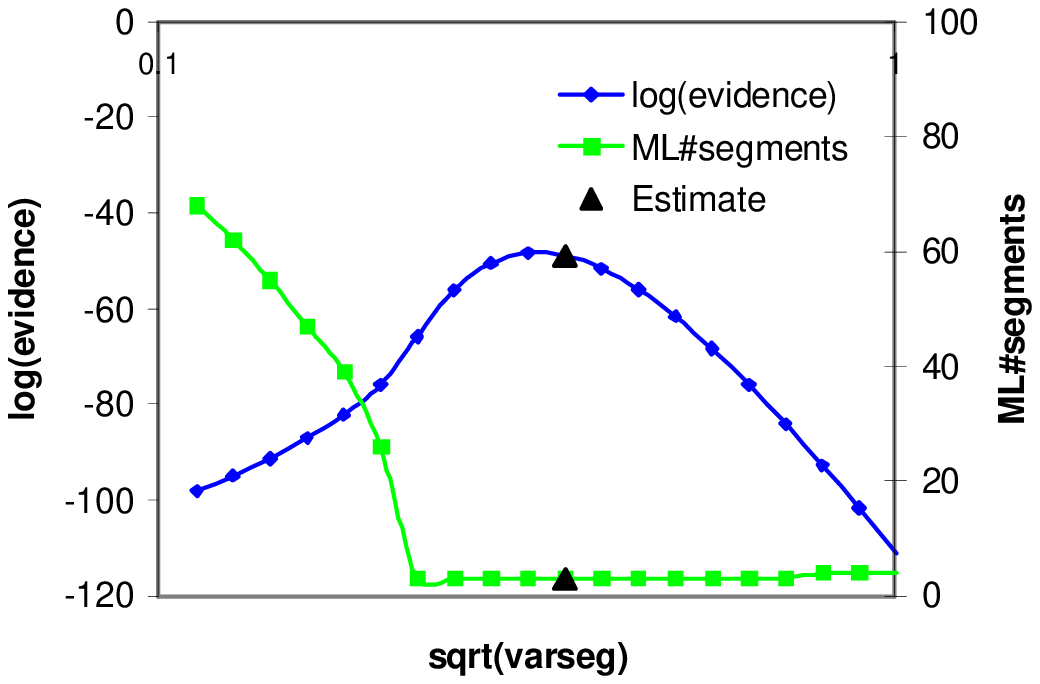}
\caption{\label{figGMFvar}[GM: medium Gaussian noise]
$\log P(\v y)$ (blue) and $\hat k$ (green) as function of $\s$ and
our estimate $\hat\s$ of $(\arg)\max_\s P(\v y)$ and $\hat
k(\hat\s)$ (black triangles).}
\end{minipage}
\hspace{0.02\textwidth}
\begin{minipage}[t]{0.49\textwidth}
\includegraphics[width=\textwidth]{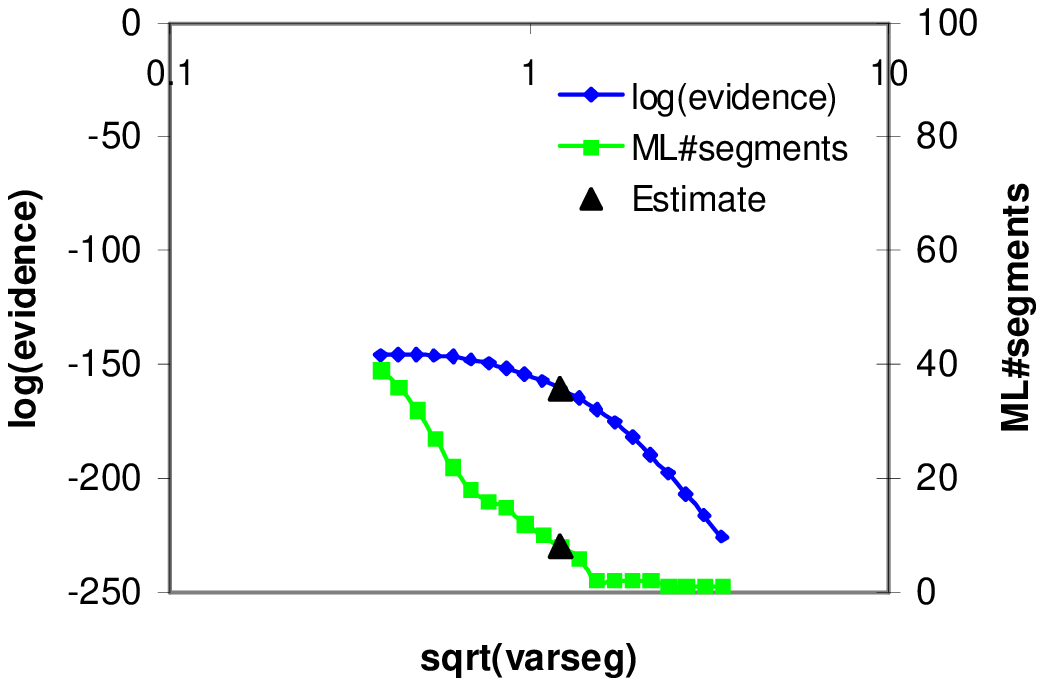}
\caption{\label{figCMFvarwG}[CMwG: medium Cauchy noise] with
Gaussian regression, $\log P(\v y)$ (blue) and $\hat k$ (green) as
function of $\s$ and our estimate $\hat\s$ of $(\arg)\max_\s P(\v
y)$ and $\hat k(\hat\s)$ (black triangles).}
\end{minipage}
\end{figure}

\paradot{Description}
In order to test our algorithm we created various synthetic data
sets. We considered piecewise constant functions with noisy
observations.
The considered function was defined $-1$ in its first quarter,
$+1$ in its second quarter, and $0$ in the last half. So the
function consists of two small and one large segments, with a
large jump at the first and a small jump at the second boundary.
For $n$ we chose 100, i.e.\ $f_1..f_{25}=-1$, $f_{26}..f_{50}=+1$,
and $f_{51}..f_{100}=0$. Data $y_t$ was obtained by adding
independent Gaussian/Cauchy noise of same scale $\s$ for all $t$.
We considered low $\s=0.1$, medium $\s=0.32$, and high $\s=1$
noise, resulting in an easy, medium, and hard regression problem
(Figures \ref{figGLPCR}-\ref{figCMFvarwG}).
We applied our regression algorithm to these 6 data sets (named
GL,GM,GH,CL,CM,CH), where we modeled noise and prior as Gaussian
or Cauchy with hyper-parameters also estimated by the Algorithms
in Table \ref{PCRegAlg}. Table \ref{tabrs} contains these and
other scalar summaries, like the evidence, likelihood, MAP segment
number $\hat k$ and their probability.

\paradot{Three segment Gaussian with low noise}
Regression for low Gaussian noise ($\s=0.1$) is very easy. Figure
\ref{figGLPCR} shows the data points $(1,y_1),..,(100,y_{100})$
together with the estimated segment boundaries and levels,
i.e.\ the Piecewise Constant Regression (PCR) curve (black).
The red curve (with the two spikes) is the posterior probability
that a boundary (break point BP) is at $t$. It is defined as
$B_t:=\sum_{p=1}^{\hat k} B_{pt}$. Our Bayesian regressor (BPCR)
is virtually sure that the boundaries are at $t_1=25$
($B_{25}=100\%$) and $t_2=50$ ($B_{25}=99.9994\%$). The segment
levels $\hat\mu_1=-0.98\approx-1$, $\hat\mu_2=0.97\approx 1$,
$\hat\mu_3=0.01\approx 0$ are determined with high accuracy i.e.\
with low deviation (green curve) $\s/\sqrt{25}=2\%$ for the first
two and $\s/\sqrt{50}\approx 1.4\%$ for the last segment. The
Bayesian regression (BR) curve $\hat\mu_t$ is identical to PCR.

\paradot{Three segment Gaussian with medium noise}
Little changes for medium Gaussian noise ($\s=0.32$). Figure
\ref{figGMPCR} shows that the number and location of boundaries is
still correctly determined, but the posterior probability of the
second boundary location (red curve) starts to get a little
broader ($B_{50}=87\%$). The regression curve in Figure
\ref{figGMBR} is still essentially piecewise constant. At $t=50$
there is a small kink and the error band gets a little wider, as
can better be seen in the (kink of the) green
$\sqrt{\Var[\mu'_t|..]}$ curve in Figure \ref{figGMPCR}.
In Figure \ref{figGMFvar} we study the sensitivity of our
regression to the noise estimate $\hat\s$. Keeping everything else
fixed, we varied $\s$ from 0.1 to 1 and plotted the log-evidence
$\log P(\v y|\s)$ and the segment number estimate $\hat k(\s)$ as
a function of $\s$. We see that our estimate $\hat\s\approx 0.35$
is close to the hyper-ML value $\s_{\text{HML}}=\arg\max_\s P(\v y|\s)
\approx 0.33$, which itself is close to the true $\s=0.32$. The
number of segments $\hat k$ is correctly recovered for a wide
range of $\s$ around $\hat\s$. If $\s$ is chosen too small (below
the critical value 0.2), BPCR cannot regard typical deviations from the
segment level as noise anymore and has to break segments into
smaller pieces for a better fit ($\hat k$ increases).
For higher noise, the critical value gets closer to $\hat\s$, but
also the estimate becomes (even) better. For lower noise, $\hat\s$
overestimates the true $\s$, but BPCR is at the same time even
less sensitive to it.

\paradot{Three segment Gaussian with high noise}
Figure \ref{figGHData} shows the data with Gaussian noise of the
same order as the jump of levels ($\s=1$). One can imagine some
up-trend in the first quarter, but one can hardly see any
segments. Nevertheless, BPCR still finds the correct boundary
number and location of the first boundary (Figure \ref{figGHPCR}).
The second boundary is one off to the left, since $y_{50}$ was
accidentally close to zero, hence got assigned to the last segment.
The (red) boundary probability curve is significantly blurred, in
particular at the smaller second jump with quite small
$B_{49}=12\%$ and $B_{50}=10\%$. The levels themselves are within
expected accuracy $\s/\sqrt{25}=20\%$ and $\s/\sqrt{50}\approx
14\%$, respectively, yielding still a PCR close to the true
function.
The Bayesian regression (and error) curve (Figure \ref{figGHBR}),
though, changed shape completely. It resembles more a local data
smoothing, following trends in the data (more on this in the next
section). The variance (green curve in Figure \ref{figGHPCR}) has
a visible bump at $t=25$, but only a broad slight elevation around
$t=50$.

\paradot{Three segment Cauchy}
The qualitative results for the Cauchy with low noise ($\s=0.1$)
are the same as for Gauss, perfect recovery of the underlying
function, and is hence not shown. Worth mentioning is that the
estimate $\hat\s$ based on quartiles is excellent(ly close to
hyper-ML) even for this low noise (and of course higher noise),
i.e.\ is very robust against the segment boundaries.

Also for medium Cauchy noise ($\s=0.32$, Figure \ref{figCMPCR})
our BPCR does not get fooled (even) by (clusters of) ``outliers''
at $t=16$, $t=48,49$, and $t=86,89,90$. The second boundary is one
off to the right, since $y_{51}$ is slightly too large. Break
probability $B_t$ (red) and variance $\Var[\mu'_t|\v y,\hat k]$
(green) are nicely peaked at $\hat t_1=25$ and $\hat t_2=51$.

For high Cauchy noise ($\s=1$, Figure \ref{figCHData}) it is nearly
impossible to see any segment (levels) at all. Amazingly, BPCR
still recovers three segments (Figure \ref{figCHPCR}), but the
first boundary is significantly displaced ($\hat t_1=14$). $B_t$
and $\Var[\mu'_t|\v y,\hat k]$ contain many peaks indicating that
BPCR was quite unsure where to break. The Bayesian regression in
Figure \ref{figCHBR} identifies an upward trend in the data
$y_{14:35}$, explaining the difficulty/impossibility of recovering
the correct location of the first boundary.

\paradot{Cauchy analyzed with Gauss and vice versa}
In order to test the robustness of BPCR under misspecification,
we analyzed the data with Cauchy noise by Gaussian BPCR (and vice
versa). Gaussian BPCR perfectly recovers the segments for low
Cauchy noise.
For medium noise (CMwG, Figure \ref{figCMPCRwG}) the outlier at
$t=49$ is not tolerated and placed in it own segment, and the last
segment is broken in two halves, but overall the distortion is
less than possibly expected (e.g.\ not all outliers are in own
segments). The reason for this robustness can be attributed to the
way we estimate $\s$. Figure \ref{figCMFvarwG} shows that the
outliers have increased $\hat\s$ far beyond the peak of $P(\v
y|\s)$, which in turn leads to lower (more reasonable) number of
segments. This is a nice stabilizing property of $\hat\s$.
The other way round, segmentation of data with medium Gaussian
noise is essentially insensitive to whether performed with
Gaussian BPCR (Fig.\ \ref{figGMPCR} and \ref{figGMBR}) or Cauchy
BPCR (GMwC, not shown), which confirms (once again) the robustness of
the Cauchy model.
But for high noise BPCR fails in both misspecification
directions.

\section{Real-World Example \& More Discussion}\label{secRWE}

\begin{figure}
\begin{minipage}[t]{0.49\textwidth}
\includegraphics[width=\textwidth]{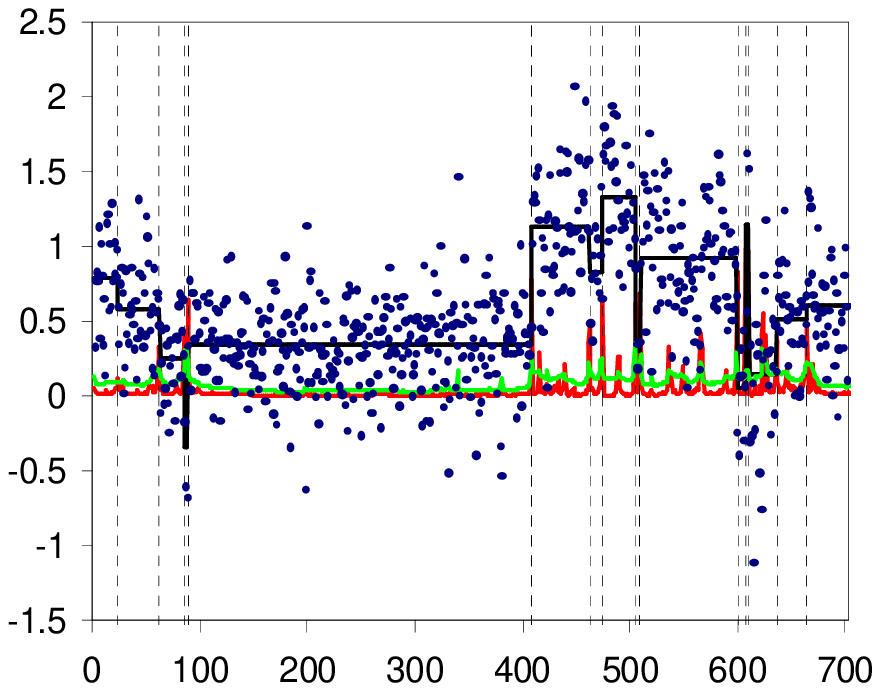}
\caption{\label{figGen31PCR}[Gen31: Aberrant gene copy \# of chromosome 1]
data (blue), PCR (black), BP (red), and variance$^{1/2}$ (green).}
\end{minipage}
\hspace{0.02\textwidth}
\begin{minipage}[t]{0.49\textwidth}
\includegraphics[width=\textwidth]{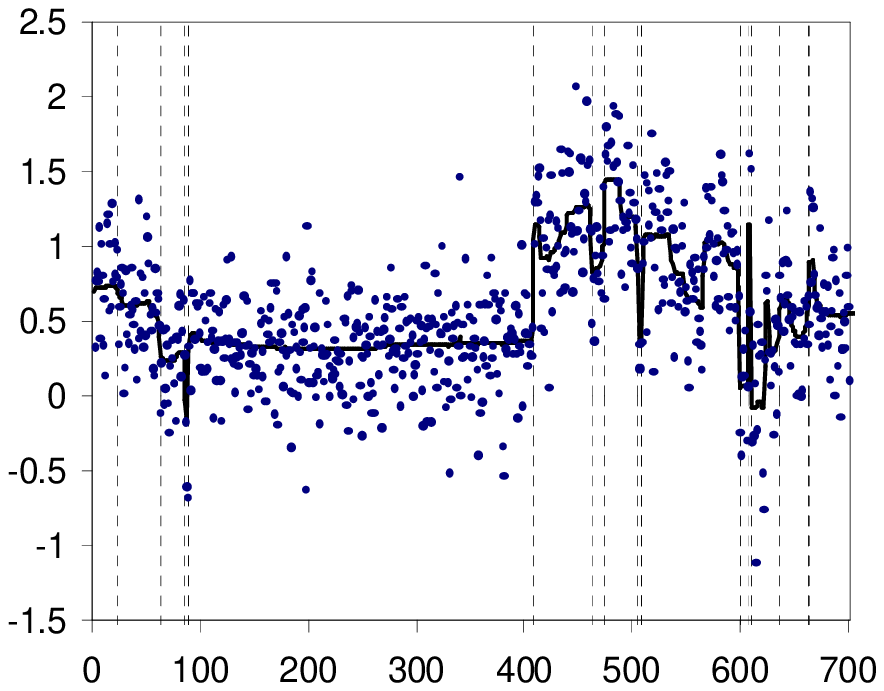}
\caption{\label{figGen31BR}[Gen31: Aberrant gene copy \# of chromosome 1]
data with Bayesian regression $\pm$ 1 std.-deviation.}
\end{minipage}
\end{figure}

\begin{figure}
\begin{minipage}[t]{0.49\textwidth}
\includegraphics[width=\textwidth]{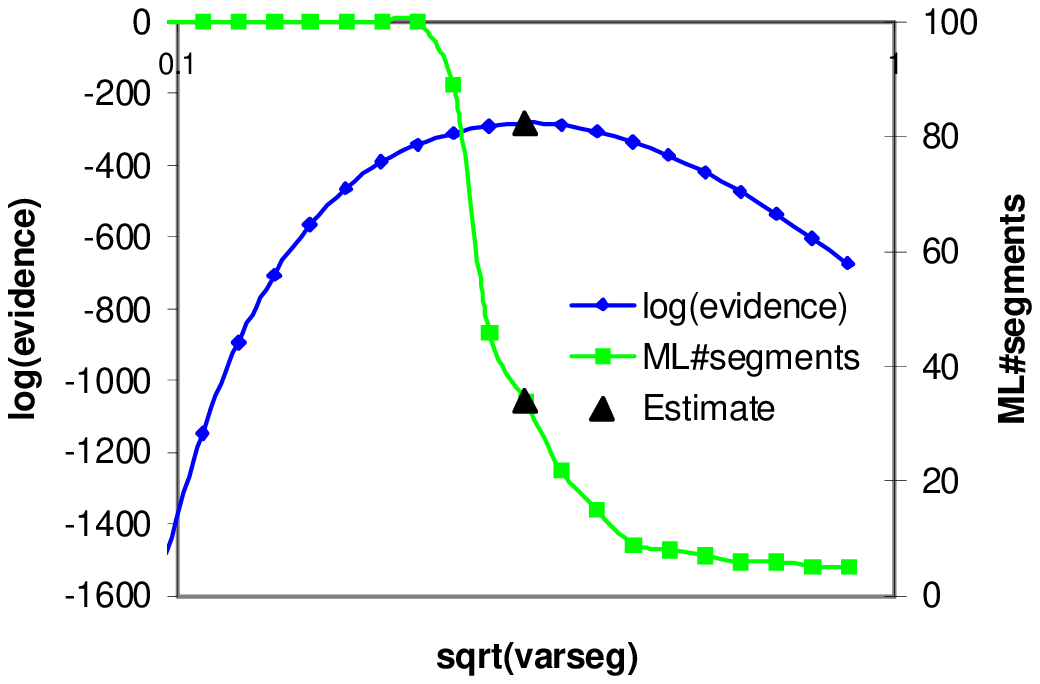}
\caption{\label{figGen31Fvar}[Gen31: Aberrant gene copy \# of chromosome 1]
$\log P(\v y)$ (blue) and $\hat k$ (green) as function of $\s$ and
our estimate $\hat\s$ of $(\arg)\max_\s P(\v y)$ and $\hat
k(\hat\s)$ (black triangles).}
\end{minipage}
\hspace{0.02\textwidth}
\begin{minipage}[t]{0.49\textwidth}
\includegraphics[width=\textwidth]{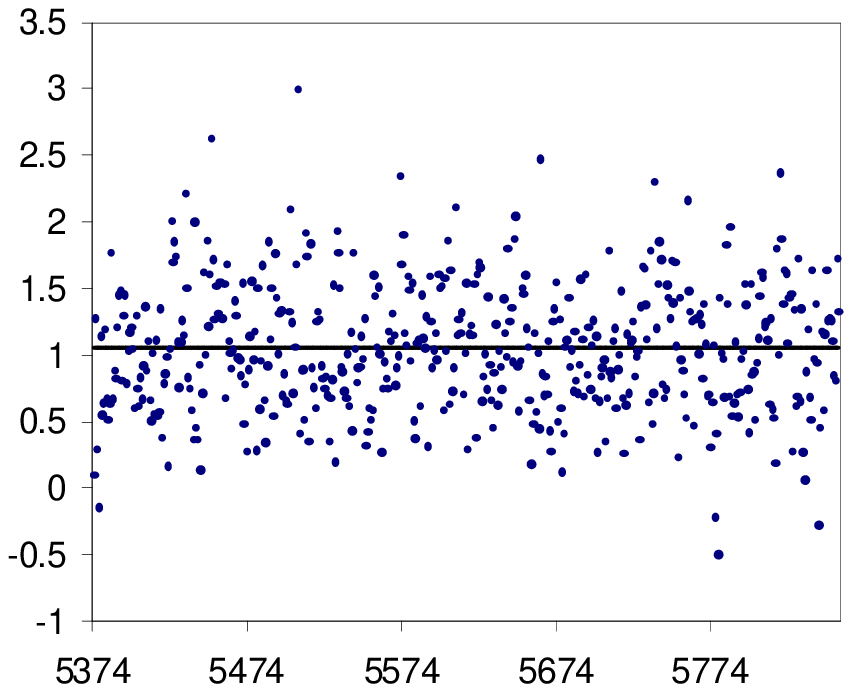}
\caption{\label{figGen59BR}[Gen59: normal gene copy \# of chromosome 9]
with Bayesian regression.}
\end{minipage}
\end{figure}

\paradot{Gene copy number data}
All chromosomes (except for the sex chromosomes in males) in a healthy
human cell come in pairs, but pieces or entire chromosomes can be lost
or multiplied in tumor cells. With modern micro-arrays one can measure
the local copy number along a chromosome. It is important to determine
the breaks, where copy-number changes. The measurements are {\em very
noisy} \cite{Pinkel:98}. Hence this is a natural application for
piecewise constant regression of noisy (one-dimensional) data. An
analysis with BPCR of chromosomal aberrations of real tumor samples,
its biological interpretation, and comparison to other methods will be
given elsewhere \cite{Hutter:06genex}. Here, we only show the
regression results of one aberrant and one healthy chromosome (without
biological interpretation).

The ``log-ratios'' $\v y$ of a normal cell (and also the
$\v\Delta$ of any cell) are very close to Gaussian distributed, so
we chose Gaussian BPCR. The log-ratios $\v y$ of chromosome 1 of a
sample known to have multiple myeloma are shown in Figure
\ref{figGen31PCR}, together with the regression results. Visually,
the segmentation is very reasonable. Long segments (e.g.\
$t=89..408$) as well as very short ones around $t=87$ and $641$ of
length 3 are detected.
The Bayesian regression curve in Figure \ref{figGen31BR} also
behaves nicely. It is very flat i.e.\ smoothes the data in long
and clear segments, wiggles in less clear segments, and has jumps
at the segment boundaries. Compare this to local smoothing
techniques \cite{Rinaldi:05}, which wiggle much more within a
segment and severely smooth boundaries. In this sense our Bayesian
regression curve is somewhere in-between local smoothing and hard
segmentation.
We also see that the regression curve has a broad dip around
$t=535..565$, although $t=510..599$ has been assigned to a single
segment. This shows that other contributions breaking the segment
have been mixed into the Bayesian regression curve. The PCR favor
for a single segment is close to ``tip over'' as can be seen from
the spikes in the break probability (red curve) in this segment.

The dependence of evidence and segment number on $\s$ is shown in
Figure \ref{figGen31Fvar}. Our estimate $\hat\s$ (black triangle)
perfectly maximizes $P(\v y|\s)$ (blue curve). It is at a deep
slope of $P(k|\v y,\s)$ (green curve), which means that
the segmentation is sensitive to a good estimate of $\hat\s$.
There is no unique (statistically) correct segmentation (number).
Various segmentations within some range are supported by
comparable evidence.

Figure \ref{figGen59BR} shows a healthy chromosome 9, correctly lumped
into one big segment.

\paradot{Posterior probability of the number of segments $P(k|\v y)$}
One of the most critical steps for good segmentation is
determining the right segment number, which we did by maximizing
$P(k|\v y)$. The whole curves shown in Figure \ref{figPKY} give
additional insight. A representative selection is presented.

For truly piecewise constant functions with $k_0\ll n$ segments
and low to medium noise, $\log P(k|\v y)$ typically raises rapidly
with $k$ till $k_0$ and thereafter decays approximately linear
(black curve). This shows that BPCR certainly does not
underestimate $k_0$ ($P(k<k_0|\v y)\approx 0$). Although it also
does not overestimate $k_0$, only $P(k\geq k_0|\v y)\approx 1$,
but $P(k_0|\v y)\not\approx 1$ due to the following reason:
If a segment is broken into two (or more) and assigned
(approximately) equal levels, the curve and hence the likelihood
does not change. BPCR does not explicitly penalize
this, only implicitly by the Bayesian averaging (Bayes factor
phenomenon \cite{Good:83,Jaynes:03,MacKay:03}). This gives very
roughly an additive term in the log-likelihood of $\odt\log n$ for
each additional degree of freedom (segment level and boundary).
This observation is the core of the Bayesian Information Criterion
(BIC) \cite{Schwarz:78,Kaass:95,Weakliem:99}.

With increasing noise, the acute maximum become more round (blue
curve), i.e.\ as expected, BPCR becomes less sure about the
correct number of segments. This uncertainty gets pronounced under
misspecification (green curve), and in particular when the true
number of segments is far from clear (or nonexistent) like in the
genome abberation example (red curve). The pink curve shows that
$\log P(k|\v y)$ is not necessarily unimodal.

\begin{table}
\begin{center}\begin{small} 
\caption[Regression summary]{\label{tabrs}
Regression summary}\vspace{1ex}
\begin{tabular}{|c|c|c|c|c|c|c|c|c|c|c|c|c|} \hline
  \rotatebox{90}{Gauss, Cauchy,} \rotatebox{90}{Low, Medium,} \rotatebox{90}{High noise, Gene}&
  \rotatebox{90}{true noise} \rotatebox{90}{scale} &
  \rotatebox{90}{data size} &
  \rotatebox{90}{method} &
  \rotatebox{90}{global mean} \rotatebox{90}{estimate} &
  \rotatebox{90}{global deviation} \rotatebox{90}{estimate} &
  \rotatebox{90}{in-segment} \rotatebox{90}{deviation est.} &
  \rotatebox{90}{log-evidence} \rotatebox{90}{$\log P(\v y)$}&
  \rotatebox{90}{rel.\ log-likelihood} \rotatebox{90}{${ll-\E[ll|\hat f]\over\Var[ll|\hat f]^{1/2}}$} &
  \rotatebox{90}{Opt.\#segm.} &
  \rotatebox{90}{Confidence} \rotatebox{90}{$P(\hat k(-1,+1)|\v y)$}
\\ \hline
Name  & $\s$ & $n$ & P &$\hat\nu$&$\hat\rho$& $\hat\s$&$\log E$&${ll-{\bf E}\over\s_{ll}}$&$\hat k$&$ C_{k(-1,+1)}$ \\ \hline\hline
GL    & 0.10 & 100 & G & -0.01 & 0.69 &  0.18 &  39 & 4.9 & $3|3$ & 74\%(0$|$20)  \\ 
GM    & 0.32 & 100 & G & -0.03 & 0.73 &  0.35 & -48 & 1.2 & $3|3$ & 44\%(0$|$29)  \\ 
GH    & 1.00 & 100 & G & -0.10 & 1.15 &  1.03 &-156 & 0.3 & $3|4$ & 13\%(10$|$12) \\ \hline
CL    & 0.10 & 100 & C & -0.02 & 0.58 &  0.09 & -17 & 1.0 & $3|3$ & 69\%(0$|$21)  \\ 
CM    & 0.32 & 100 & C & -0.09 & 0.70 &  0.27 &-127 & 0.8 & $3|3$ & 38\%(0$|$27)  \\ 
CH    & 1.00 & 100 & C & -0.20 & 0.99 &  0.86 &-234 & 0.9 & $3|4$ & 12\%(11$|$11) \\ \hline
GMwC  & 0.32 & 100 & C &  0.00 & 0.49 &  0.17 & -70 & 1.5 & $3|3$ & 27\%(0$|$26)  \\ 
CMwG  & 0.32 & 100 & G &  0.01 & 1.24 &  1.22 &-160 & 2.9 & $5|8$ &  8\%(8$|$8)   \\ \hline
Gen31 &  --  & 769 & G &  0.55 & 0.45 &  0.30 &-283 &-1.5 &$15|34$&  6\%(6$|$6)   \\ 
Gen59 &  --  & 483 & G &  1.05 & 0.47 &  0.44 &-336 &-2.3 & $1|1$ &  8\%(0$|$6)   \\ \hline
\end{tabular}
\end{small}\end{center}
\end{table}

\paradot{Miscellaneous}
Table \ref{tabrs} summarizes the most important quantities of the
considered examples.

While using the variance of
$\v\Delta$ as estimate for $\hat\s$ tends to overestimate $\s$ for
low noise, the quartile method does not suffer from this
(non)problem.

The usefulness of quoting the evidence cannot be overestimated.
While the absolute number itself is hard to comprehend,
comparisons (based on this absolute(!) number) are invaluable.
Consider, for instance, the three segment medium Gaussian noise
data $y_{\text{GM}}$ from Figure \ref{figGMPCR}. Table \ref{tabrs}
shows that $\log E($GM$)=-48$, while $\log E($GMwC$)=-70$, i.e.\
the odds that $\v y_{\text{GM}}$ has Cauchy rather than Gaussian
noise is tiny $\e^{48-70}< 10^{-9}$, and similarly the odds that
$\v y_{\text{CM}}$ has Gaussian rather than Cauchy noise is
$\e^{127-160}<10^{-14}$. This can be used to decide on the model
to use. For instance it clearly indicates that noise in Gene31 and
Gen59 is not Cauchy for which log-evidences would be $-398$ and
$-406$, respectively. The smallness of the relative
log-likelihoods does not indicate any gross misspecification.

The indicated 4$^{th}$ segment for GH and CH is spurious, since it
has length zero (two breaks at the same position). In Gene31, only
15 out of the indicated 34 segments are real. The spurious ones
would be real had we estimated the breaks $\v{\hat t}$ jointly,
rather than the marginals $t_p$ separately. They would often be
single data segments at the current boundaries, since it costs
only a single extra break to cut off an ``outlier'' at a boundary
versus two breaks in the middle of a segment.

In the last column we indicated the confidence $C_{\hat k}$
$(C_{\hat k-1},C_{\hat k+1})$ of BPCR in the
estimate $\hat k$. For clean data (GL,GM,CL,GM) it is certain that
there are at least 3 segments. We already explained the general
tendency to also believe in higher number of segments.

\section{Extensions \& Outlook}\label{secMisc}

The core Regression($\v A,n,k_{max}$) algorithm does not care
where the in-segment evidence matrix and moments $\v A$ come from.
This allows for plenty of easy extensions of the basic idea.

If the segment levels are known to belong to a discrete set (e.g.\
integer DNA copy numbers \cite{Picard:05segcl}), this simply
corresponds to a discrete prior on $\mu$ and leads naturally to a
Grid sum (rather than by need) as in EstGeneral().

If each segment can have its own (unknown) variance $\s_m^2$,
we can assume some prior over $\s_m$ and average \req{Adef} (which
depends on $\s_m$, notationally suppressed) additionally over
$\s_m$. Possibly $P(\s_m|...)$ depends on some hyper-parameter that
now has to be estimated instead of $\s$; all the better if not.

We assumed a constant regression function within a segment.
Actually any other function could be used. We simply choose
likelihood and prior for a single segment and compute its evidence
$A_{ij}^0$. This is all what Regression() needs to determine the
segment number and boundaries. Once we have the segment boundaries
it is easy to compute the in-segment quantities we are
interested in, e.g.\ the MAP or mean regression curve.

For instance, if we consider all linear functions within a segment,
we get a piecewise linear regression curve. But note that this curve
is not continuous. This model is, for instance good, if the true
function is essentially piecewise constant, but there is an additional
underlying trend (slope) in the segments. Using non-linear functions allows
to handle more complicated trends.

Piecewise linear (or other) {\em continuous} regression is more
complicated. Assume that $\mu_p$ in \req{Qrec} does not denote the
level of the whole segment $p$, but its level at the right
boundary, which together with $\mu_{p-1}$ determines the linear
function in segment $p$. Only after fixing $\mu_p$, left and right
side decouple. So the recursion analogous to \req{Qrrec} now
involves a quantity $Q$ which in addition to $(i,j)$ also depends
on $(\mu_l,\mu_m)$. This functional recursion may approximately be
solved by discretizing $\{(\mu_l,\mu_m)\in\SetR^2\}$, or by
approximating $Q$ by a 2-dimensional Gaussian in $(\mu_l,\mu_m)$
and storing only the 2 means and the $2\times 2$ covariance matrix
for each $(i,j)$.
%
The following two simpler heuristic approaches may work
sufficiently well in practice: One could ignore the continuity
constraint when determining the boundaries, and only take them
into account in the subsequent (much simpler) regression problem
with known boundaries. Another possibility is to consider instead
of the continuous piecewise linear function $f$ its piecewise
constant derivative $f'$, i.e.\ use BPCR on $\Delta_t$ and
finally integrate the result.

It is also not necessary to use a parametric model for the
noise. If different segments can have different noise
distributions, we could compute the in-segment evidence, mean, and
variance $A_{ij}^r$ based on some (fast) non-parametric model. If
all segments have the same distribution, we could
non-parametrically estimate a single density for the differences
$\v\Delta$ and then deconvolve the density (e.g.\ by
$\mbox{FFT}^{-1}(\sqrt{\mbox{FFT(density)$\!\!$}}\;$), and
henceforth use this as prior for $\s$ in EstGeneral(). As
non-parametric density estimator we could use the fast (linear-time) exact
Bayesian tree model \cite{Hutter:04bayestree}.

Finally, for (very) large $n$, say $>1000$, the $O(k_{max}n^2)$
algorithm is too slow. Fortunately, there is nearly no interaction between
distant segments; boundary $t_k$ is often practically
independent of where $t_{k\pm 2}$, $t_{k\pm 3}$, etc.\ are placed.
This suggests to break the whole data set into smaller overlapping
pieces, where each piece should be long enough to contain at least
four segments. Then boundaries $t_2^{piece},...,t_{k-2}^{piece}$
of each piece are used, and appropriately merged. For the Bayesian
regression curve one should use some blending on the
overlap. If single segments are very long, one could coarsen
(locally lump together) the data and later refine around the
boundaries.

\section{Summary}\label{secDisc}

We considered Bayesian regression of piecewise constant functions
with unknown segment number, location and level.
We derived an efficient algorithm that works for any noise and
segment level prior, e.g.\ Cauchy which can handle outliers. We
derived simple but good estimates for the in-segment variance. We
also proposed a Bayesian regression curve as a better way of
smoothing data without blurring boundaries. The Bayesian approach
also allowed us to straightforwardly determine the global
evidence, break probabilities and error estimates, useful for
model selection and significance and robustness studies. We
discussed the performance on synthetic and real-world examples.
Many possible extensions have been discussed.

\paradot{Acknowledgements}
Thanks to IOSI for providing the gene copy \# data and to Ivo
Kwee for discussions.


\begin{small}

\end{small}
\end{document}